\documentclass[a4paper, 11pt]{article}
\usepackage[english]{babel}
\usepackage{amsmath, amssymb, amsfonts, bm}
\usepackage{mathrsfs}
\usepackage{textcomp}
\usepackage[mathcal]{euscript}
\usepackage{enumerate}
\usepackage{mathbbol}
\usepackage{float}
\usepackage{color,graphicx,epstopdf}
\usepackage{dsfont}
\usepackage{hyperref}
\usepackage[noadjust]{cite}
\usepackage{caption}
\usepackage{subcaption}
\usepackage{mathtools}
\usepackage{tikz-cd}
\usepackage{mathtools}
\usepackage{algorithm}
\usepackage{algpseudocode}
\usepackage{url}
\usepackage{authblk}

\newtheorem{definition}{Definition}

\newtheorem{lemma}{Lemma}

\newtheorem{proposition}{Proposition}

\newcommand{\bfu}{{\mathbf u}}
\newcommand{\bfU}{{\mathbf U}}
\newcommand{\bfx}{{\mathbf x}}
\newcommand{\bfX}{{\mathbf X}}

\newcommand{\bff}{{\mathbf f}}
\newcommand{\bpi}{{\bm \pi}}
\newcommand{\mcalV}{{\mathcal V}}
\newcommand{\mcalT}{{\mathcal T}}
\newcommand{\mcalS}{\mathcal S}
\newcommand{\mcalU}{{\mathcal U}}
\newcommand{\mcalX}{{\mathcal X}}

\usepackage[left=1.5cm,right=1.5cm,top=2cm,bottom=2cm]{geometry}
\usepackage{float}

\title{\bf How Cooperation and Competition Arise in Regional Climate Policies: RICE as a Dynamic Game}         
\author[1]{Yijun Chen}
\author[1]{Guodong Shi} 

\affil[1]{Australian Center for Field Robotics, School of Aerospace, Mechanical and Mechatronic Engineering, The University of Sydney, NSW, Australia, email: \{yijun.chen, guodong.shi\}@sydney.edu.au}
\date{}
\begin{document}
	\maketitle
	
	\begin{abstract} 
	One of the most widely used models for studying the geographical economics of climate change is the Regional Integrated model of Climate and the Economy (RICE). In this paper, we investigate how cooperation and competition arise in regional climate policies under the RICE framework from the standpoints of game theory and optimal control. First, we show that the RICE model is inherently a dynamic game. Second, we study both cooperative and non-cooperative solutions to this RICE dynamic game. In cooperative settings, we investigate the global social welfare equilibrium that maximizes the weighted and cumulative social welfare across regions. We next divide the regions into two clusters: developed and developing, and look at the social welfare frontier under the notion of Pareto optimality.  We also present a receding horizon approach to approximate the global social welfare equilibrium for robustness and computational efficiency. For non-cooperative settings, we study best-response dynamics and open-loop Nash equilibrium of the RICE game. A Recursive Best-response Algorithm for Dynamic Games (RBA-DG) is proposed to describe the sequences of best-response decisions for dynamic games, which indicates convergence to open-loop Nash equilibrium when applied to the RICE game by numerical studies. We also study online receding horizon feedback decisions of the RICE game. A Receding Horizon Feedback Algorithm for Dynamic Games (RHFA-DG) is proposed. All these proposed solution concepts are implemented and open sourced using the latest updated parameters and data. The results reveal how game theory may be used to facilitate international negotiations towards consensus on regional climate-change mitigation policies, as well as how cooperative and competitive regional relations shape climate change for our future.
	\end{abstract}

	\section{Introduction}
	The issue of global warming has emerged as a central international environmental question over decades. As a consequence of industrialization and economic development, human-caused emissions of greenhouse gases, most notably carbon dioxide ($\rm CO_2$), contribute to a significant increase in the mean atmospheric temperature by $1.15^\circ$C relative to the pre-industrial age \cite{portner2022climate}. This temperature deviation yields significant changes in the global climate and ecosystem, including increasing larger wildfires \cite{barbero2015climate},  sea level rise \cite{courchamp2014climate}, melting of ice lands \cite{salim2021review}, and so on. To assess the damages of anthropogenic greenhouse gases, especially $\rm CO_2$, on society, scientists in this field employ the Integrated Assessment Models (IAMs) \cite{hope2013critical,anthoff2013uncertainty,nordhaus2013dice}. IAMs simulate the dynamics of the economy-climate interactions by incorporating mathematical models from both economics and geophysical science. The Dynamic Integrated model of Climate and Economy (DICE) \cite{nordhaus2013dice,nordhaus2017revisiting}, developed by William Nordhaus, who received the $2018$ Nobel Memorial Prize in Economic Sciences in large part for this body of work, is one of the most well-known IAMs. The DICE model is a globally aggregated model and treats global warming as a single-agent problem. Considering the crucial aspect of regional socio-economic heterogeneity, Nordhaus also proposed the Regional Integrated model of Climate and the Economy (RICE) \cite{nordhaus1996regional,nordhaus2010economic}, which is a decentralized version of the DICE model. By dividing the world into several regions, the RICE model takes the vantage point in determining how multiple regions may jointly design climate policies and cope with the global warming issue.
	
	 Global society has been making efforts in developing sensible strategies to achieve international consensus on climate-change mitigation over the last several decades. There are several important international agreements on climate change. In $1992$, the United Nations Framework Convention on Climate Change (UNFCC) was the first formal global treaty to be signed by $154$ (now $198$) parties to explicitly address climate change, which established an annual forum for international discussions aimed at stabilizing greenhouse gas concentrations in the atmosphere \cite{bodansky1993united}. These meetings produced the well-known Kyoto Protocol \cite{bohringer2003kyoto}, Copenhagen Accord \cite{rogelj2010copenhagen}, and Paris Agreement \cite{schleussner2016science}. The Kyoto Protocol was a binding agreement that was signed in $1997$ and expired in $2012$. The Kyoto Protocol called for developed countries to reduce emissions by an average of $5$ percent below pre-industrial levels, and established a system to monitor countries’ progress \cite{bohringer2003kyoto}. The meeting in Copenhagen in $2009$ was called to establish a replacement for the Kyoto Protocol. Although it failed to establish binding emission limits after $2012$, countries recognized ``the scientiﬁc view that the increase in global temperature should be below $2^\circ$C'' \cite{rogelj2010copenhagen}. The Paris Agreement was signed in $2016$ and called for all countries to set emissions-reduction pledges/targets, with the goals of preventing the mean atmospheric temperature from rising $2^\circ$C above pre-industrial levels and pursuing efforts to keep it below $1.5^\circ$C \cite{schleussner2016science}. Despite intensified diplomacy in these meetings, most existing climate-change treaties are neither sufficient nor mandatory, some of which even stalled due to a lack of political will \cite{nwanguma2016analysis,rhodes2017us}. Due to economic competition and political divide, an international enforceable agreement on specific emission-reduction control has not yet been reached.
	
	 There are two coupled but conflicting sides in regional emission reduction policies: Regions are affected by the same global climate system; they also decide their climate strategies to benefit their individual economic benefits and political self-interests. Therefore, it is actually a decision-making process where competition occurs. Game theory has been a fundamental tool in explaining decision-making of independent and self-interested players in a strategic and competitive setting \cite{von2007theory}. In a game, each player takes its action to maximize its payoff function. The payoff function of each player relies on not only its own action, but also other players’ actions, resulting in inherent competitions. A common solution concept in game theory is Nash equilibrium under which no one will benefit by changing its action when others remain unchanged \cite{holt2004nash}. In the RICE model, the regional climate strategies are associated with a common global climate dynamic, and as a result, the climate policy decisions fall into the concept of dynamic games established from the interface between game theory and optimal control theory \cite{bacsar1998dynamic}. In a dynamic game, the strategic interaction among players recurs over time. The group of players is associated with a dynamical game state that depends on all players' actions. The goal is for each player to take an action to maximize each player's cumulative payoff function over time that depends on the game state, its own action, and other players' actions.
	
	In this paper, we investigate how cooperation and competition in climate policy across the various regions in the globe under the RICE framework affect the formation of international climate treaties, the implementation of regional climate-change mitigation measures, and the resulting implications of competitive/non-cooperative decisions for climate change. In control community, there are a few efforts on studying the climate-change mitigation measures. The work of \cite{kellett2019feedback} provided a tutorial introduction to the DICE model and proposed a receding horizon approach to DICE. A bi-objective optimal control problem (OCP) on DICE was studied, objectives of which are maximizing social welfare and minimizing atmospheric temperature deviation \cite{heris2020multi}. The work of  \cite{carlino2020multi} studied a multi-objective stochastic OCP on DICE, which accounts for stochastic disturbances and aligns with physical targets posed by international agreements on climate change mitigation. 
	
	The contributions of this work are summarized as follows. We show that the RICE model is inherently a dynamic game, termed the RICE game. Both cooperative and non-cooperative solutions to this RICE game are considered:
	\begin{itemize}
		\item	For cooperative solutions, we study global social welfare maximization problem where all regions take actions that maximize the weighted and cumulative social welfare across all regions, which serves as a benchmark. Next, we classify regions into two clusters of developed regions and developing regions, and consider the concept of Pareto equilibrium, describing that any attempt to benefit one cluster by deviating to some outcome will necessarily result in a loss in satisfaction of the other cluster. Our simulation results show that under different Pareto equilibria, the social welfare of developed regions and developing regions would not drastically change, and the atmospheric temperature deviation is quite robust. Finally, we apply a receding horizon approach to approximate the solution to the global social welfare maximization problem. The receding horizon control presents a favorable approximation property. 
		
		\item	For non-cooperative solutions, we study best-response dynamics and open-loop Nash equilibrium of the RICE game. Multiple plays/episodes of dynamic games are considered, where each player chooses the sequence of control for the next episode that maximizes its social welfare given other players' sequences of control in the current episode.A Recursive Best-response Algorithm for Dynamic Games (RBA-DG) is proposed. By applying it to the RICE game, simulation results show that the obtained sequence of actions converges to a steady point, indicating that RBA-DG is useful for computing the open-loop Nash equilibrium of the RICE game. We also study online receding horizon feedback decisions of the RICE game. A single play of dynamic games is considered, where players observe the current game state and other players' actions, and then apply a receding-horizon feedback decision-making approach to predict the action in the future. A Receding Horizon Feedback Algorithm for Dynamic Games (RHFA-DG) is proposed. By implementing it over the RICE game, simulation results show that with the nature of competition, receding horizon, and myopic assumption about other players' actions, the emission-reduction rates become lower in most time steps. 
	\end{itemize}
	
	The implementations of our studies are based on the RICE-2010 model with a few parameters updated from the latest updated data. There have been a few implementations in the literature on DICE and RICE \cite{nordhaus2020website,faulwasser2018mpc,faulwasser2018github,kellett2016github,anthoff2021github}. As for the DICE model, Nordhaus provided a GAMS implementation of the DICE-$2013R$ model with the objective of solving a social welfare maximization problem \cite{nordhaus2013dice,nordhaus2020website}. Faulwasser, Kellett and Weller published a Matlab implementation of the DICE-$2013R$ model, which replicated the basic functionality of GAMS implementation \cite{kellett2016github}. Further, Faulwasser, Kellett and Weller released a Matlab implementation of the MPC-DICE model which uses model predictive control to approximately solve the DICE OCP. As for the RICE model, Nordhaus provided an Excel implementation (without optimization module), and a GAMS implementation (which is currently inaccessible) \cite{nordhaus2020website}. Anthoff and Errickson created an implementation in Julia programming language \cite{anthoff2021github}, which re-coded the Excel version of RICE without optimization module \cite{anthoff2021github}. Our implementations of various solution concepts for RICE under the proposed dynamic-game perspective are developed relying on these previous efforts. We have open-sourced our implementation as a RICE-GAME framework, with a Matlab and Casadi-based implementation of RICE dynamic game, Preprint at \url{arXiv.org}, code for download at \cite{chen2022matlab}.

	The paper is organized as follows. Section~\ref{sec:preliminaries} provides preliminaries of the DICE/RICE model and dynamic games. Section~\ref{sec:RICE-GAME}
	represents the RICE model as a dynamic game. Section~\ref{sec:cooperative-solutions} considers cooperative solutions for the RICE game under three cooperative settings: global social welfare maximization, RICE Pareto frontier, and receding-horizon global social welfare maximization. Two non-cooperative settings are then considered. Best-response dynamics and open-loop Nash equilibrium for the RICE game are studied in Section~\ref{sec:BR-NE}. A receding horizon feedback planning approach for RICE game is proposed in Section~\ref{sec:RHFD}. The paper ends with concluding remarks in Section~\ref{sec:conclusions}.
	
	\section{Preliminaries}\label{sec:preliminaries}
	
	In this section, we introduce some preliminary knowledge on the DICE/RICE model and dynamic games. 
	
	\subsection{The DICE Model}\label{subsec:DICE}
	
\begin{figure}[h]
		\centering
		\includegraphics[width = 0.5\linewidth]{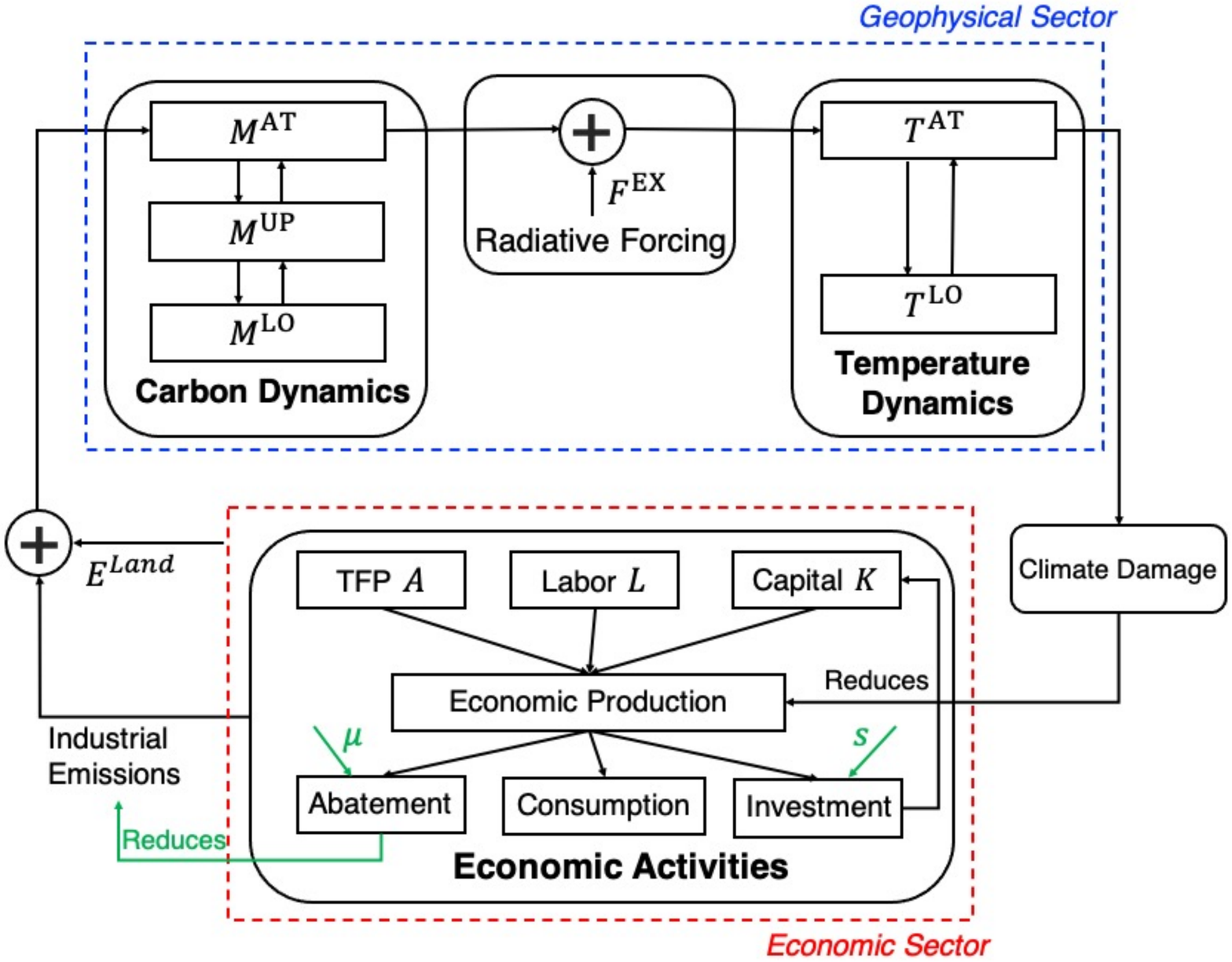}
		\caption{A block diagram for the DICE model.}
		\label{fig:workflow_dice}
	\end{figure}

	The DICE model \cite{nordhaus2013dice} is one of integrated
	assessment models (IAMs) that simulate the interplay between economy and climate, and quantify the social cost of $ \rm CO_2$ emissions. The DICE model operates in periods of $5$ years and its latest version starts from the year $2015$ as the initial year. The DICE model is composed of two sectors (see  Fig.~\ref{fig:workflow_dice}): a geophysical sector (blue dotted block) that accounts for the global interaction between carbon  and temperature, and an economic sector that is globally aggregated for the world total (red dotted block). 
	
	{\noindent \bf Geophysical and Economic Sectors}. In the geophysical sector, the DICE model considers  ${\rm CO}_{2}$ emissions as the major contributor to climate change. The geophysical sector is constructed as follows. 
	\begin{itemize}
		\item  There are two main sources of $\rm CO_{2}$ emissions: industrial $\rm CO_{2}$ emissions related to the carbon intensity (denoted by $\sigma$) of global economic activities and natural $\rm CO_{2}$ emissions due to land use changes, $E^{\rm land}$.  The global $\rm CO_{2}$ emissions as the sum of industrial and natural emissions drive the carbon cycle of the Earth. 
		\item The Carbon dynamics are described by  a three-reservoir model \cite{post1990global} on the  carbon flows among three reservoirs: the atmosphere, the upper oceans and  biosphere, and the deep oceans. The average carbon masses in those reservoirs are represented by $M^{\rm AT}$, $M^{\rm UP}$, and $M^{\rm LO}$, respectively. 
		\item Accumulations of $\rm CO_{2}$ emissions and other greenhouse gases warm the Earth's surface through enhanced radiative forcing.  Radiative forcing resulted by $\rm CO_{2}$ emissions has a logarithmic dependence on the atmospheric carbon mass; greenhouse gases other than $\rm CO_{2}$ emissions contribute to exogenous radiative forcing $F^{\rm EX}$.
		\item  The rise in temperature at the Earth's surface is driven by radiative forcing. Temperature dynamics are captured by a two-layer model \cite{schneider1981atmospheric}. Given the temperature in year $1750$ as zero reference, $T^{\rm AT}$ and $T^{\rm LO}$ represent the temperature deviation in the atmosphere and in the lower ocean from those of the reference year, respectively. 
	\end{itemize}
	The economic sector of DICE is based on the Cobb-Douglas production function \cite{cobb1928theory}, where gross economic output is determined by total factor productivity $A$, labor $L$, and capital $K$.  The total factor productivity and labor evolve exogenously; the capital dynamics follows the Solow-Swan model \cite{swan1956economic} where capital depreciates over time and is replenished by investment. 
	
	\medskip
	
	{\noindent \bf Climate-Economy Feedback}. The DICE model establishes two feedback loops between the geophysical sector and the economic sector: (i) the industrial $ {\rm CO}_2$ emissions are a by-product of economic activities; (ii)  the atmospheric temperature rise  has an negative impact on economic production. 
	
	\medskip
	
	{\noindent \bf Control Inputs}.  The DICE model assumes two control decisions: the saving rate $s$ and the emission-reduction rate $\mu$. The saving rate $s$ represents  the ratio of investment to the economic output; the emission-reduction rate $\mu$   represents the rate at which industrial $\rm CO_{2}$ emissions are reduced. By adjusting the saving rate, it is possible to balance consumption today and consumption in the future. By increasing the emission-reduction rate to slow down $ {\rm CO}_2 $ emissions as a ``climate investment'',  the  currently available amount of consumption and investment will be reduced \cite{nordhaus2013dice}. This climate investment will lower climate damage and therefore potentially  increase consumption in the future. 
	
	\medskip
	
	{\noindent \bf System Outputs}. The economic output is counted as output net of emission abatement cost and climate damage. Social welfare is calculated as the discounted sum of the population-weighted utility of per capita consumption.
	

	\medskip
	
	{\noindent \em DICE Variables.} All variables described above are time-dependent, although not explicitly written. The variables $T^{\rm AT}, T^{\rm LO}, M^{\rm AT},  M^{\rm UP}, M^{\rm LO}$ and $F^{\rm EX}$ belong to the geophysical sector, whereas the variables $K,A, L, \sigma,$ and $ E^{\rm land}$ belong to the economic sector. Some variables  evolve independently, whereas others evolve in an interdependent manner. The variables that evolve independently as exogenous signals are $F^{\rm EX}, A, L, \sigma,$ and $ E^{\rm land}$; the variables evolving in an interdependent manner are $T^{\rm AT}, T^{\rm LO}, M^{\rm AT},  M^{\rm UP}, M^{\rm LO}, $ and $K.$
	
	\subsection{The RICE Model}\label{subsec:RICE}
	
	The RICE model is a variant of the DICE model that accounts for regional climate damages and control decisions \cite{nordhaus2010economic,nordhaus2011estimates}. Since being proposed in 1990s,  calibration of the  RICE model has been updated several times, and the latest version of the RICE model is the RICE-2011 model on which our study is based. The RICE-2011 model uses a time-step of $10$ years, starting from the year 2005 as the initial year. 
	
\begin{figure}
		\centering
		\includegraphics[width = 0.5\linewidth]{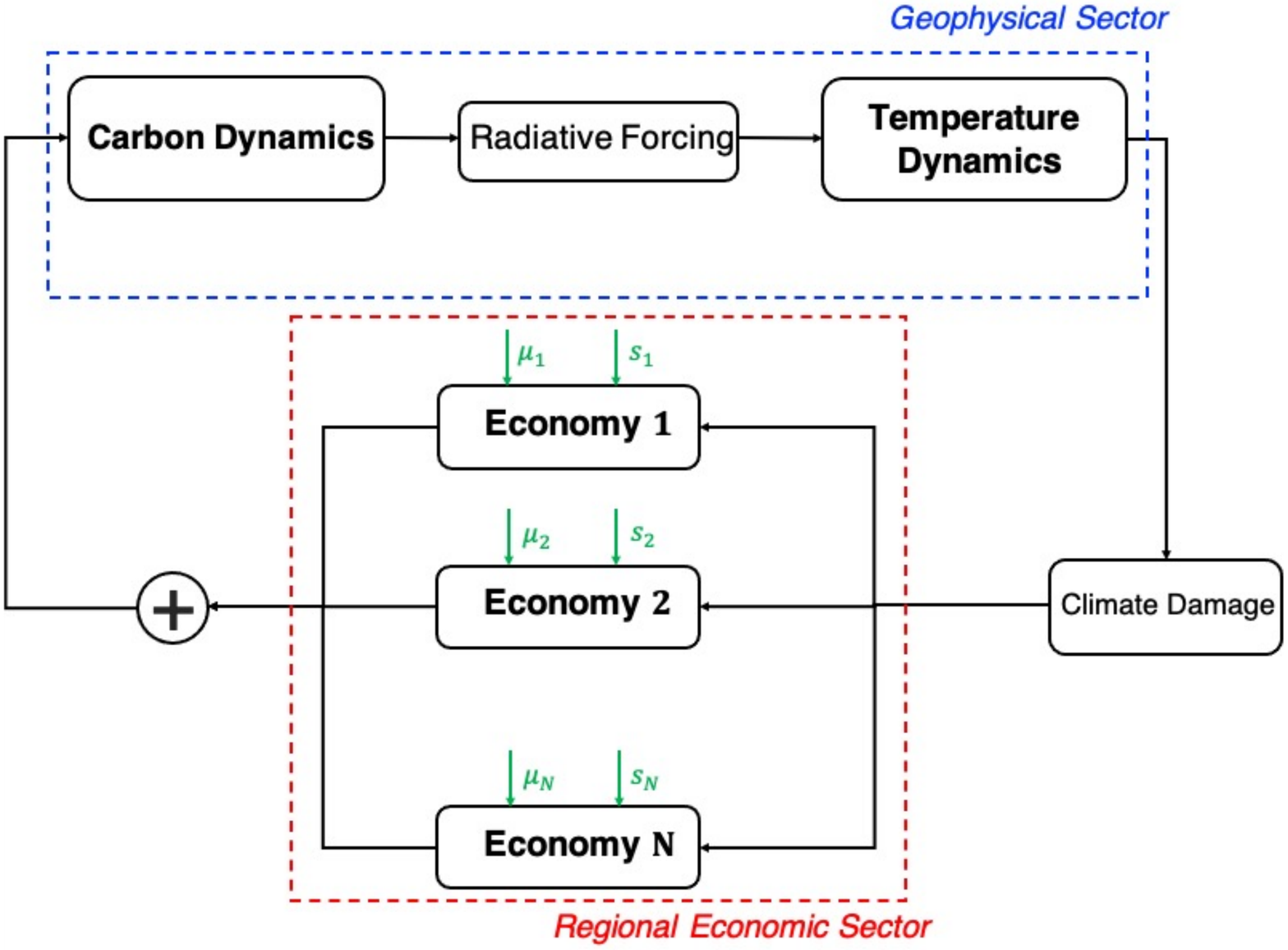}
		\caption{The workflow of the RICE-2011 model.}
		\label{fig:bd_rice}
	\end{figure}

	The RICE-2011 model integrates a global geophysical sector with regional economic sectors (see Fig.~\ref{fig:bd_rice}):
	\begin{itemize}
		\item The global geophysical sector of the RICE-2011 model contains the same carbon dynamics and temperature dynamics as the DICE model;
		\item The regional economic sectors of the RICE-2011 model disaggregates the world into $12$ regions (US, EU, Japan, Russia, Non-Russian Eurasia, China, India, Middle East, Africa, Latin America, other high income countries, and other Asian countries), each of which is equipped with region-specific  climate damage level, economic factors, and saving rate and emission-reduction rate as local control inputs. 
	\end{itemize}
	\medskip

	{\noindent \it RICE Variables.} The variables $T^{\rm AT}, T^{\rm LO}, M^{\rm AT},  M^{\rm UP},$ $M^{\rm LO},$ and $F^{\rm EX}$ in the global geophysical sector of the RICE-2011 model are inherited from the DICE model, while the variables $K_{i},A_{i}, L_{i}, \sigma_{i},$ and $ E^{\rm land}_{i}, i\in\{1,2,\dots,12\},$ in the regional economic sectors of the RICE-2011 model correspondingly depend on specific  regions. 
	
	\subsection{Dynamic Games}\label{subsec:DTDG}
	
	The theory of dynamic games lies in the interface between   game theory and optimal control, which involves a dynamic decision process for multiple players \cite{bacsar1998dynamic}. An $n$-player discrete-time dynamic game over a finite horizon is defined as follows. 
	
	\medskip
	
	{\noindent \bf Dynamic Game.} The $n$ players are indexed in $\mcalV:= \{1,2,\dots,n\}$; time is discrete with the steps  indexed in $\mcalT:= \{0,1,\dots,T\}$. Each player can manipulate the game through its control decisions, and the control decision space of player $i \in \mcalV$ is denoted by $\mcalU_{i} \subseteq \mathbb{R}^{d}$. At each time step $t=0,\dots,T$, the   decision executed by player $i$  is denoted  by $\bfu_{i}(t) \in \mcalU_{i}$. We also use $\bfu(t) = [\bfu_{1}^{\top}(t);\dots;\bfu_{n}^{\top}(t)]$, $\bfU_{i} = [\bfu_{i}^{\top}(0);\dots;\bfu_{i}^{\top}(T)]$ and $\bfU = [\bfU_{1};\dots;\bfU_{n}]$ to represent the all-player decision profile at time $t$, the player-$i$ decision throughout the time horizon, and the decision profile for all players and for all time stpes. The control decisions of all players excluding player $i$ at time step $t$ is denoted by $\bfu_{-i}(t)$, and the control decisions of all players excluding player $i$ over the entire horizon is represented  by $\bfU_{-i}$.
	
	For each $t=0,1,\dots,T$, the group of players are associated with a dynamical state    $\bfx(t)\subseteq \mathbb{R}^{m}$ that evolves according to
	\begin{equation}\label{eq:DTDG_dynamics}
		\bfx(t+1) = \bff(t, \bfx(t),\bfu(t)), \quad \bfx(0) = \bfx_{0}, \quad t=0,\dots,T,
	\end{equation}
	with $ \bfx_{0}$ being the initial state. At each time $t=0,\dots, T$, upon playing $\bfu_{i}(t)$, the agent $i$   receives a payoff $g_{i}(\bfx(t),\bfu_{i}(t),\bfu_{-i}(t)) \in \mathbb{R}$ given other players' actions $\bfu_{-i}(t)$ and the current state $\bfx(t)$, where  $g_{i}(\bfx(t),\bfu_{i}(t),\bfu_{-i}(t))$ is a continuous function with respect to $\bfx(t)$, $\bfu_{i}(t)$, and $\bfu_{-i}(t)$. The cumulative payoff of agent $i$ throughout the time horizon is therefore 
	\begin{equation}\label{eq:DTDG_payoff_function}
		\mathsf{J}_i(\bfX,\bfU_{i},\bfU_{-i}) = \sum_{t=0}^{T}g_{i}(\bfx(t),\bfu_{i}(t),\bfu_{-i}(t))
	\end{equation}
	where $\bfX=(\bfx^\top(0), \dots, \bfx^\top(T))^\top$. 
	Each player's goal is to make decisions for maximizing its cumulative payoff function; the system dynamics produces a terminal state $\bfx(T+1)$ towards the end of the time horizon as a result of those decisions. 
	
	In what follows, we present cooperative and non-cooperative settings as well as three solution concepts for discrete-time dynamic game: Pareto solution, open loop Nash equilibrium (NE) and feedback Nash equilibrium.
	
	\subsubsection{Cooperative Decisions: Pareto Solutions}
	In the cooperative setting,  players are able to communicate and cooperate with each other to achieve their objectives, and all  players know  the system dynamics and payoff functions of other players. Pareto optimality is an efficiency concept \cite{mas1995microeconomic} under which any attempt to benefit one player by deviating to some other outcome will necessarily result in a loss in satisfaction of another player. 
	\begin{definition}[Pareto Efficiency]\label{df:pareto}
		A   decision profile $\bfU^{p}$ is Pareto efficient for the dynamic game if there does not exist another decision profile  $\bfU$ such that  
		\begin{itemize}
			\item[(i)] there holds for all $i \in \mathcal{V}$ $
			\mathsf{J}_i(\bfX^{p},\bfU_{i}^{p},\bfU_{-i}^{p})\leq
			\mathsf{J}_i(\bfX,\bfU_{i},\bfU_{-i})$;
			\item[(ii)] there exists at least one $k\in \mathcal{V}$ such that $
			\mathsf{J}_{k}(\bfX^{p},\bfU_{k}^{p},\bfU_{-k}^{p})<
			\mathsf{J}_{k}(\bfX,\bfU_{k},\bfU_{-k})$
		\end{itemize}		
		Here $\bfX$ and $\bfX^{p}$ are the states evolved under $\bfU$ and $\bfU^{p}$, respectively. The set of all Pareto solution is called the Pareto frontier.
	\end{definition}
	
	The following Lemma provides a convenient way of computing Pareto solutions \cite[Lemma 6.1]{engwerda2005lq}.
	
	\begin{lemma}\label{lemma:pareto}
		Consider a set of parameters
		$\mathcal{C}:=\{\mathbf{c} = (c_{1},\dots,c_{n}):\ c_{i}\geq 0, \text{\rm and } \sum_{i=1}^{n}c_{i}=1\}.$
		If a   decision profile $\bfU^{p}$ 
		is such that 
		\begin{equation}\label{eq:pareto}
			\bfU^{p} \in \arg\max_{\bfU}\Big\{ \sum_{i=1}^{n} c_{i}\mathsf{J}_i\Big\},
		\end{equation}
		for some $\mathbf{c} \in \mathcal{C}$, then $\bfU^{p}$ is Pareto efficient.
	\end{lemma}

	\subsubsection{Non-Cooperative Setting: Nash Equilibriums}
	In the non-cooperative setting of games, Nash equilibriums (NE) mark one of the most important solution concepts. For a dynamic game, the information structure in terms of what players know before a decision is made at a particular time becomes critical in properly defining NE. There are two basic types of information structure: the open loop information structure and the feedback information structure. 
	
	In the open loop information structure, each player knows  the initial state $\bfx_{0}$, and then plans at $t=0$ all the control decisions      $\bfu_{i}(t)$ for $t \in \mcalT$. Consequently, the open loop control decision of  $\bfu_{i}(t)$ can be written as 
	\begin{equation}\label{eq:ol_control}
		\bfu_{i}(t) = \bfu_{i}(t,\bfx_{0}).
	\end{equation} 
	Denote an open loop  decision profile by $\bfU^{\ast}$ where $\bfU_{i}^{\ast} := [\bfu^{\ast}_{i}(0,\bfx_{0});\dots;\bfu^{\ast}_{i}(T,\bfx_{0})], i \in \mcalV$. We introduce the following definition. Here with slight abuse of notation we also write  $\mathsf{J}_i(\bfX,\bfU_{i},\bfU_{-i})$ as  $\mathsf{J}_i(\bfx_{0},\bfU_{i},\bfU_{-i})$ noting the fact that $\bfX$ is uniquely determined by $\bfx_{0}$ and $(\bfU_{i},\bfU_{-i})$. 
	\begin{definition}[Open Loop NE]\label{def:ol_ne}
		Given the initial state $\bfx_{0}$, a control decision profile $\bfU^{\ast}$ is said to be an open loop Nash equilibrium control decision profile if there holds for all $i \in \mcalV$ and all $ \bfU_{i}$ that
		\begin{equation}\label{eq:ol_ne}
			\mathsf{J}_i(\bfx_{0},\bfU_{i}^{\ast},\bfU_{-i}^{\ast}) \geq \mathsf{J}_i(\bfx_{0},\bfU_{i},\bfU_{-i}^{\ast}).
		\end{equation}
	\end{definition}
	
	In the feedback information structure, at time step $t \in \mcalT$, each player knows the current system state $\bfx(t)$. As a consequence, each player $i \in \mcalV$ may employ a feedback law $\bpi_{i}:\mcalT \times \mcalX \to \mcalU_{i}$ to determine its decisions by
	\begin{equation}\label{eq:fb_control}
		\bfu_{i}(t) = \bpi_{i}(t,\bfx(t)).
	\end{equation}
	In this case, we denote  the overall feedback law as  $ \bpi  = [\bpi_{1};\dots;\bpi_{n}]$, and then write $\mathsf{J}_i(\bfX,\bfU_{i},\bfU_{-i})$ as  $\mathsf{J}_i(\bfx_{0},\bpi_{i},\bpi_{-i})$ noting $(\bfX,\bfU_{i},\bfU_{-i})$ is uniquely determined by $\bfx_{0}$ and $\bpi$. We introduce the following definition. 
	\begin{definition}[Feedback NE]\label{def:fb_ne}
		For any initial state $\bfx_{0}$,  a feedback profile $ \bpi^{\ast}$ is said to be a  feedback Nash equilibrium   for the dynamic game if there holds for all $i\in \mathcal{V}$ and all $\bpi_{i}$ that
		\begin{equation}\label{eq:fb_ne}
			\mathsf{J}_i(\bfx_{0},\bpi^{\ast}_{i},\bpi^{\ast}_{-i}) \geq \mathsf{J}_i(\bfx_{0},\bpi_{i},\bpi^{\ast}_{-i}).
		\end{equation}
	\end{definition}

	\section{RICE as a Dynamic Game}\label{sec:RICE-GAME}
	
	In this section, we show that the RICE model is inherently a dynamic game where regional saving rates and emission-reduction rates regulate global temperature, and then the global temperature has impact on regional social welfares through climate damage. Our presentation is based on the RICE-2011 model with slight modifications, but the nature of being a dynamic game is embedded in all RICE models.

	There are $12$ regions in the   RICE-2011 model. Each region is considered a player and the regions are  indexed in $\mcalV = \{1,2,\dots,n\}$ with $n=12$. We operate the RICE dynamic game in periods of $5$ years, starting from the year 2020 as the initial year\footnote{The  RICE-2011 model operates in periods of $10$ years starting from 2005.}. Taking  the discrete time step index $\mcalT = \{0,1,\dots,T\}$, the relation between an actual calendar year and the corresponding discrete time step is determined by  
	\begin{equation}\label{eq:time_step}
		year(t) = year(0) + 5t, \quad year(0) =2020.
	\end{equation}
	Note that although most variables in the RICE-2011 model are defined as flows per year and only some variables are in flows per decade \cite[Supplementary Material]{nordhaus2011estimates}, all variables in this present of the RICE dynamic game are defined as flows per year. 
	
	\subsection{System Dynamics}
	We define the dynamical  state   of the RICE model  at time step $t \in \mcalT$ as
	\begin{align} \label{eq:rice_state_variables}
		\bfx(t) = [T^{\rm AT}(t); T^{\rm LO}(t); M^{\rm AT}(t);  M^{\rm UP}(t); M^{\rm LO}(t); K_{1}(t); \dots; K_{n}(t)] \in \mathbb{R}^{n+5}.& 
	\end{align}
	Let the  control decision of region $i \in \mcalV$ at time step $t \in \mcalT$ be
	\begin{equation}\label{eq:rice_control_i}
		\bfu_{i}(t) = [s_{i}(t);\mu_{i}(t)]^\top:=[\bfu_{i[1]}(t);\bfu_{i[2]}(t)]^\top \in [0,1]^{2}.
	\end{equation} 
	Consequently, the control decisions of the RICE dynamic game at time step $t \in \mcalT$ of all players are 
	\begin{equation}\label{eq:rice_control_whole}
		\bfu(t) = [s_{1}(t);\mu_{1}(t);\dots;s_{n}(t);\mu_{n}(t)] \in [0,1]^{24}.
	\end{equation} 
	According to the RICE model, the dynamics of $\bfx(t)$ can be written as 
	\begin{equation}\label{eq:dynamics}
		\bfx(t+1) = \bff(t, \bfx(t),\bfu(t)), \quad \bfx(0) = \bfx_{0}, \quad t \in \mcalT,
	\end{equation}
	where $\bff:=[f_{1};f_{2};\dots;f_{n+5}]^\top$ follows from the interdependency among the geophysical signals and the economic signals, and the feedback between the two sectors. In what follows, we briefly  describe the form of entries in the dynamics $\bff$. For a more detailed description, please refer to \cite{nordhaus2011estimates}.

	\medskip

	{\noindent \bf Carbon dynamics.} There are three carbon reservoirs: the atmosphere, the upper oceans and the biosphere, and the deep oceans. The atmospheric carbon reservoir has an additional input, the global ${\rm CO}_{2}$ emissions $E(t)$ that is related to economic activities and land use at time $t$.   The carbon dynamics for carbon transition among the three reservoirs  are described by
	\begin{equation}\label{eq:carbon_dynamics}
		\begin{bmatrix}
			M^{\rm AT}(t+1) \\
			M^{\rm UP}(t+1) \\
			M^{\rm LO}(t+1)
		\end{bmatrix}
		=
		\begin{bmatrix}
			\zeta_{11}&\zeta_{12}&0\\
			\zeta_{21}&\zeta_{22}&\zeta_{23}\\
			0&\zeta_{32}&\zeta_{33}\\
		\end{bmatrix}
		\begin{bmatrix}
			M^{\rm AT}(t) \\
			M^{\rm UP}(t) \\
			M^{\rm LO}(t)
		\end{bmatrix}
		+
		\begin{bmatrix}
			\xi_{1}\\0\\0
		\end{bmatrix} E(t),
	\end{equation}
	where  $\zeta_{11},\zeta_{12},\zeta_{21},\zeta_{22},\zeta_{23},\zeta_{32},\zeta_{33}$ and $\xi_{1}$ are constant parameters. 
	
	\medskip

	{\noindent \bf Temperature dynamics.} The evolution of the atmospheric and ocean temperature is governed by  the following equations:
	\begin{equation}\label{eq:temperature_dynamics}
		\begin{bmatrix}
			T^{\rm AT}(t+1) \\
			T^{\rm LO}(t+1)
		\end{bmatrix}
		=
		\begin{bmatrix}
			\phi_{11}&\phi_{12}\\
			\phi_{21}&\phi_{22}
		\end{bmatrix}
		\begin{bmatrix}
			T^{\rm AT}(t) \\
			T^{\rm LO}(t)
		\end{bmatrix}
		+
		\begin{bmatrix}
			\xi_{2}\\0
		\end{bmatrix} F(t),
	\end{equation}
	where  $\phi_{11},\phi_{12},\phi_{21},\phi_{22}$ and $\xi_{2}$ are constant parameters, and the  radiative forcing  at time step $t$, $F(t)$, is computed as
	\begin{equation}\label{eq:radiative_forcing}
		F(t)=\eta \log_{2}(\frac{M^{ \rm AT}(t)}{M^{\rm AT,1750}}) + F^{\rm EX}(t). 
	\end{equation}	
	Here $\eta$ and $M^{\rm AT,1750}$ are constants, and $ F^{\rm EX}(t)$ represents radiative forcing due to
	other greenhouse gases at time step $t$.

	\medskip
	
	{\noindent \bf Economic dynamics.} The economy of each region $i \in \mcalV$ at time step $t \in \mcalT$ follows from the   Cobb-Douglas production function \cite{cobb1928theory}:
	\begin{equation}\label{eq:gross}
		Y_{i}(t) = A_{i}(t)K_{i}(t)^{\gamma_{i}}(t)L_{i}(t)^{1-\gamma_{i}}, 
	\end{equation}
	where  $Y_{i}(t),A_{i}(t),K_{i}(t)$ and $L_{i}(t)$ represent region $i$'s gross economic output, total factor productivity, capital stock and labor at time step $t$, respectively. Here, $\gamma_{i}, i \in \mcalV,$ are constant parameters. 
	
	\medskip
	
	{\noindent \bf Economy-climate feedback}.  Global ${\rm CO}_{2}$ emissions at time step $t$ are the sum of natural emissions because of each region's land use at time step $t$, $E^{\rm land}_{i}(t)$, and industrial emissions resulted from each region's economic activities at time step $t$. Each region $i$'s industrial emissions depend on each region $i$'s carbon intensity at time step $t$, $\sigma_{i}(t)$, which is an exogenous variable. Consequently, global ${\rm CO}_{2}$ emissions at time step $t$ is described by
	\begin{equation}\label{eq:total_emission}
		E(t) = \sum_{i=1}^{n} \Big(\sigma_{i}(t)(1-\mu_{i}(t))Y_{i}(t) +E^{\rm land}_{i}(t) \Big),
	\end{equation}
	where $\mu_{i}(t), i \in \mcalV, t \in \mcalT,$ are control decisions representing the emission-reduction rate.

	The emission abatement cost fraction as the percentage of gross economic output spent on emission-reduction effort
	at time step $t$ is given by 
	\begin{equation}\label{eq:abt_cost}
		\Lambda_{i}(t) =  1-\theta_{i}^{[1]}(t) \mu_{i}(t)^{\theta_{i}^{[2]}},
	\end{equation}
	where $\theta_{i}^{[2]}, $ and $\theta_{i}^{[1]}(t), i \in \mcalV, t \in \mcalT,$ are parameters. 
		The parameters $\theta_{i}^{[1]}(t), i \in \mcalV, t \in \mcalT,$ are calculated by 
		\begin{equation}\label{eq:theta_1}
			\theta_{i}^{[1]}(t) = \frac{pb_{i}}{1000\cdot\theta_{i}^{[2]}}(1-\delta^{pb}_{i})^{t-1}\cdot\sigma_{i}(t),
		\end{equation}
		where $pb_{i}, i \in \mcalV,$ represent the price of backstop technology at time step $t = 0$ for region $i$ to replace all carbon fuels, and $\delta^{pb}_{i}, i \in \mcalV,$ are constant parameters. 
	The damage function 	$\Omega_{i}(t)$ is the percentage of gross economic output damaged by temperature rising. As a result, the net economic  output  $Q_{i}(t)$ (economic output after the emission-reduction spending and climate damage), is given by
	\begin{equation}\label{eq:net}
		Q_{i}(t) = \Omega_{i}(t) \Lambda_{i}(t) Y_{i}(t).
	\end{equation}
	The Solow-Swan model \cite{swan1956economic}  gives a description of capital accumulation of each region $i \in \mcalV$:
	\begin{equation}\label{eq:captial_dynamics}
		K_{i}(t+1)=(1-\delta^{\rm K}_{i})^{5} K_{i}(t)+5 s_{i}(t) Q_{i}(t),
	\end{equation}
	where $\delta^{\rm K}_{i}, i \in \mcalV,$ are constant parameters, and $s_{i}(t), i \in \mcalV, t \in \mcalT,$ are control decisions representing the saving rate, i.e., the percentage of net economic output invested in capital.

	\subsection{Damage Functions}
	The rising atmospheric temperature has a negative impact on economic production. Although there are various specifications offering various estimates of the damage function, there are no substantial discrepancies among them \cite{botzen2012sensitive}. In the RICE-2011 model, damages to economic gross output caused by rising temperatures are considered to be region-specific and dependent on factors such as atmospheric temperature deviation, sea level rises, and atmospheric carbon mass. However, the form of the damage function of the RICE-2011 model is not given explicitly. In the current study, a simplified form of the damage function is employed, which only depends on the rising atmospheric temperature deviation:  
	\begin{equation}\label{eq:dmg_fr}
		\Omega_{i}(t) = 1 - a_{i}^{[1]}T^{\rm AT}(t) - a_{i}^{[2]}T^{\rm AT}(t)^{a_{i}^{[3]}},
	\end{equation}
	where $a_{i}^{[1]},a_{i}^{[2]},$ and $a_{i}^{[3]},i\in \mcalV,$ are parameters calibrated to yield a certain amount of damage loss to regional economic gross output. Each region's damage loss to its regional economic gross output at $2^\circ C$ is presented in Table~\ref{tb:table_damage_loss}. For example, India would suffer a damage loss of $1.55\%$  to its economic gross output when the atmospheric temperature deviation reaches  $2^\circ C$.
	
	\begin{table}[htbp]
		\centering
		\small
		\setlength\tabcolsep{2pt}
		\begin{tabular}{|c|c|c|c|c|c|c|}
			\hline
			& US & EU & Japan & Russia & Eurasia & China\\ 
			\hline
			Loss & 0.56\% & 0.64\%&0.65\%&0.46\%&0.52\%&0.66\%\\
			\hline
			& India & MidEast & Africa & LatAm & OHI & OthAsia\\
			\hline
			Loss &1.55\%&1.19\%&1.47\%&0.66\%&0.62\%&1.04\%\\
			\hline
		\end{tabular}	
		\caption{Each region's damage loss with respect to its gross economic output when atmospheric deviation is at $2^\circ C.$ }
		\label{tb:table_damage_loss}	
	\end{table}

	\subsection{Payoff Functions}
	
	Note that from the net economic output $Q_i(t)$,  a total amount of $s_i(t) Q_i(t)$ has been made as  investment. The remaining part $C_i(t)=(1-s_i(t))Q_i(t)$ can then be used for consumption. In RICE models, for region $i$, the social welfare of the population $L_i(t)$ consuming $C_i(t)$ of economic output at time $t$ is defined by the population-weighted utility of per capita consumption 
	\begin{equation}\label{eq:utility}
		g_{i}(C_{i}(t),L_{i}(t)) =  L_{i}(t) \cdot \frac{(\frac{C_{i}(t)}{L_{i}(t)}) ^{1-\alpha_{i}}-1}{1-\alpha_{i}},
	\end{equation}
	where $\alpha_{i}$ is a constant. The cumulative social welfare of  region $i$ across the time horizon   is then given by
	\begin{align}\label{eq:social_welfare}
		\mathsf{J}_i &=  \sum_{t=0}^{T} \frac{g_{i}(C_{i}(t),L_{i}(t))}{(1+\rho_{i})^{5t}} \nonumber\\
		&= \sum_{t = 0}^{T} \Bigg( \Big[ \frac{A_{i}(t)L_{i}(t)^{1+\alpha_{i}-\gamma_{i}}}{(1-\alpha_{i})(1+\rho_{i})^{5t}}\Big(1-\bfu_{i[2]}(t)\Big)\Big(1 - a_{i}^{[1]}\bfx_{1}(t) - a_{i}^{[2]}\bfx_{1}(t)^{a_{i}^{[3]}}\Big)\Big(1-\theta_{i}^{[1]}(t) \bfu_{i[1]}(t)^{\theta_{i}^{[2]}}\Big)\bfx_{5+i}(t)^{\gamma_{i}}  \Big]  
\nonumber\\
		& \ \ \ - \frac{L_{i}(t)}{(1-\alpha_{i})(1+\rho_{i})^{5t}} \Bigg)\nonumber\\
		&:=\mathsf{J}_i(\bfX,\bfU_{i},\bfU_{-i})
	\end{align}
	where $\rho_{i}$ is a constant discounting factor. For each region  $i \in \mcalV$, naturally it will attempt to maximize its cumulative social welfare. 
	
	We have now formally represented the RICE-2011 model as a dynamic game, termed the {\it RICE game} where  the regions as players seek to  plan their control decisions in emission-reduction rate and saving rate for the entire time horizon (600 years) $\bfU_{i} = [s_{i}(0);\mu_{i}(0);\dots;s_{i}(T);\mu_{i}(T)]$ so as to maximize their payoff functions~\eqref{eq:social_welfare}  subject to the underlying dynamical system \eqref{eq:dynamics}, represented in \eqref{eq:temperature_dynamics}-\eqref{eq:captial_dynamics}. In what follows when we implement the proposed RICE game, the values of the initial state $\bfx(0)$ and the parameters  are updated in the following way: the initial state are calibrated to match the data in the year $2020$; the parameters in the geophysical sector use the latest updated values in the DICE-2016 model \cite{nordhaus2017revisiting}, while the parameters in the regional economic sector remain unchanged.

	\subsection{The Social Cost of $\rm CO_2$}
	The social cost of  $\rm CO_2$ (SCC) is a central concept for understanding and implementing climate change policies.
	In a market setting like a cap-and-trade regime, the SCC would serve as the trading price of carbon emission permits. In a carbon-tax regime, the SCC would be the carbon tax for the regions who want to emit carbon emissions. The SCC in a particular year is defined as the decrease in aggregate consumption in that year that would change the current expected value of social welfare by the same amount as a one unit increase in carbon emissions in that year \cite{newbold2013rapid}.  The regional SCC is then given by 
	\begin{align}
		{\rm SCC}_{i}(t)&= -1000 \cdot \frac{\partial \mathsf{J}_i}{\partial E_{i}(t)}\Big/ \frac{\partial \mathsf{J}_i}{\partial C_{i}(t)} \notag \\ 
		&= -1000 \cdot \frac{\partial C_{i}(t)}{\partial E_{i}(t)} 	\label{eq:scc}.
	\end{align}

	\section{The RICE Game: Cooperative Solutions}\label{sec:cooperative-solutions}
	
	In this section, we study the solutions to the RICE game under cooperative settings. First of all, we revisit the classical RICE solution concept  defined by  a system-level    social welfare maximization, and present numerical  solutions under the current calibration of climate damage function and the updated geophysical and economic parameters and conditions. Next, we move to Pareto solutions to the proposed RICE game, and present the Pareto frontier of the RICE game between developing and developed regions. Finally, we introduce a receding horizon solution to the classical  RICE solution as the counter-part for MPC-DICE developed in \cite{faulwasser2018mpc,faulwasser2018github}.

	\subsection{RICE Social Welfare Maximization}\label{subsec:RICE-SWM}
	
	The RICE-2011 model focused on the sum of the weighed regional social welfare across all regions:
	\begin{equation}\label{eq:weighed_social_welfare_over_regions}
		\mathsf{W}_{\mathbf{c}} =  \sum_{i = 1}^{n} c_{i} \mathsf{J}_i,
	\end{equation}
	where  $0<c_{i}<1$ is know as  the Negishi weight for region $i$ with $\sum_{i=1}^n=1$.
		The values of the $c_i$ were calibrated in the work of \cite{nordhaus2011estimates} to be the inverse of the marginal utilities of consumption when maximizing~\eqref{eq:weighed_social_welfare_over_regions} subject to RICE dynamics~\eqref{eq:carbon_dynamics}-\eqref{eq:captial_dynamics} under the assumption that there is no abatement of ${\rm CO_2}$ emissions \cite{nordhaus1996regional,nordhaus2013dice,stanton2011negishi}.

\subsubsection{Solution Concept}		One benchmark cooperative solution to the RICE game is for a centralized climate policy planner to compute the $\bfU_{1},\dots,\bfU_{n}$ for all regions that achieve the maximal value of $\mathsf{W}_{\mathbf{c}}$, for a given initial condition $\mathbf{x}_0$.

	\begin{definition}\label{def:swm}
		A decision profile  $\bfU^{\rm w}$ is a global social welfare equilibrium if it is a solution to the following optimization problem
		\begin{equation}\label{eq:swm}
			\begin{aligned}
				&  \max_{\bfU_{1},\cdots,\bfU_{n}}
				& &  \mathsf{W}_{\mathbf{c}} (\bfX, \bfU) \\
				&  {\rm subject \ to}
				& & \bfx(t+1) = \bff(t, \bfx(t),\bfu(t)), \ \  \bfx(0) = \bfx_{0},  t \in \mcalT \\
				&&& \bfu(t)   \in [0,1]^{24}, \ \  t \in \mcalT. 
			\end{aligned}
		\end{equation}
	\end{definition}

	\subsubsection{Results}We run a total of  $120$ $5-$year periods, and therefore we set $T = 120$. We obtain the global social welfare equilibrium  $\bfU^{\rm w}$ for the RICE game under Definition \ref{def:swm} by solving the corresponding OCP~\eqref{eq:swm}.

\begin{figure}
		\centering
		\begin{subfigure}{0.65\textwidth}
			\includegraphics[width=\textwidth]{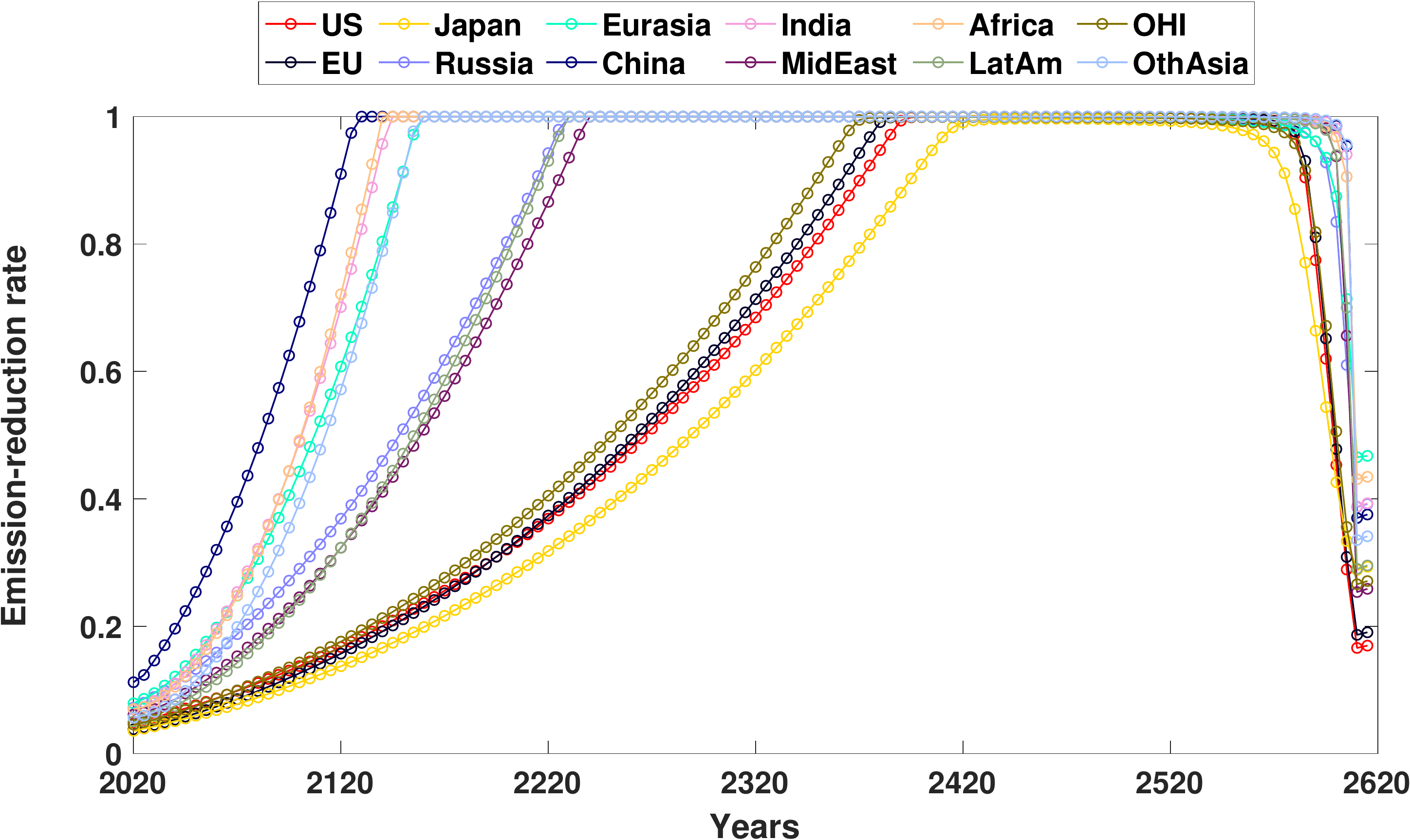}
			\caption{Optimal emission-reduction rate.}
			\label{fig:cop_reduction_rate}
		\end{subfigure}
		\begin{subfigure}{0.65\textwidth}
			\includegraphics[width=\textwidth]{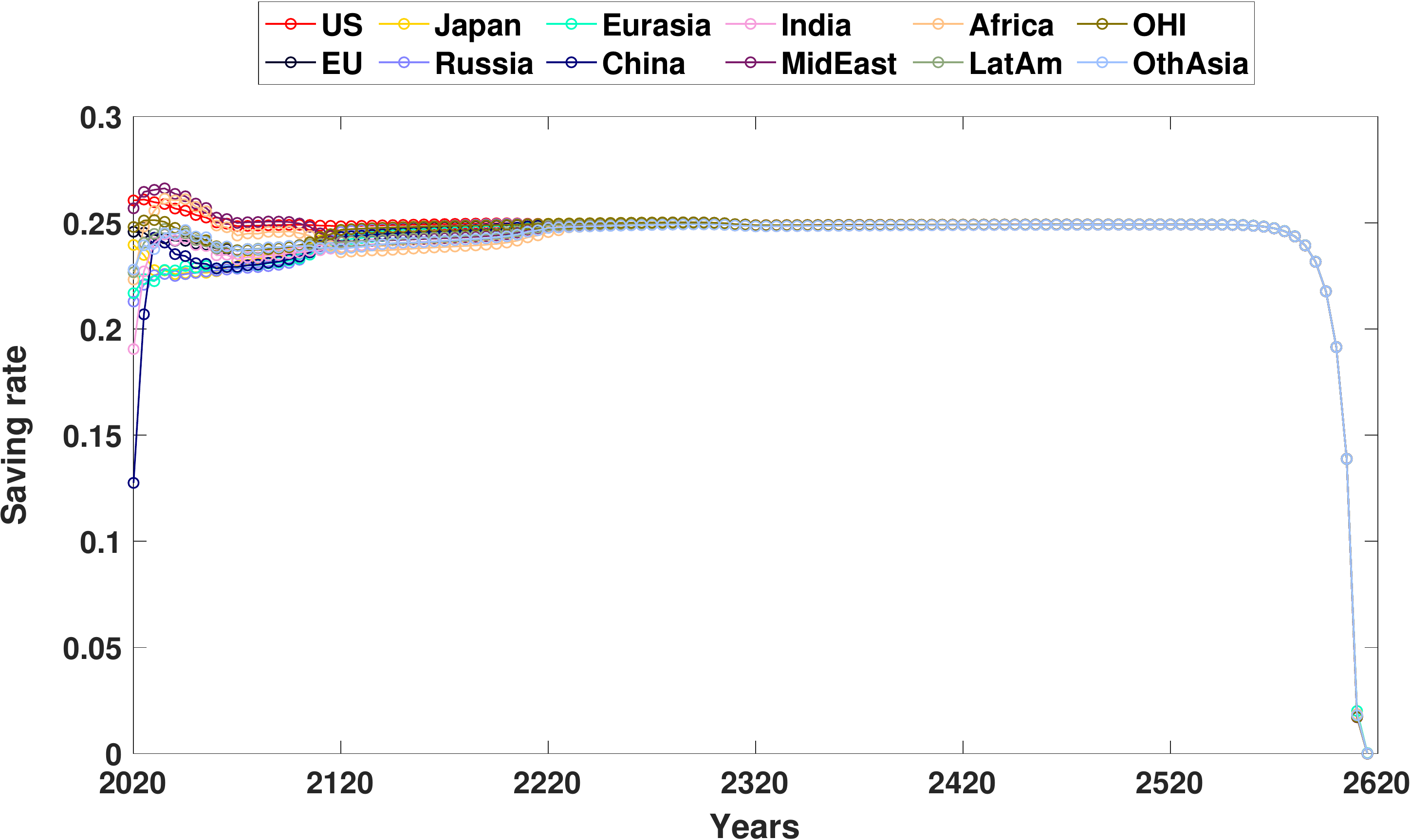}
			\caption{Optimal saving rate.}
			\label{fig:cop_saving_rate}
		\end{subfigure}
		\caption{Each region's optimal emission-reduction rate and saving rate under $\bfU^{\rm w}$.}
		\label{fig:dg_cd}
	\end{figure}

First,	each region's optimal control decisions in   emission-reduction rate and saving rate under $\bfU^{\rm w}$ are plotted  in Fig~\ref{fig:dg_cd}.  Clearly,  regions including  China, Africa, India, Eurasia, OthAsia, and Russia have relatively high emission-reduction rates who should reach carbon neutral status (i.e., the emission-reduction rate $\mu_i(t)$ becoming exactly $1$) relatively sooner. The reason for that could be twofolds. First, regions such as Africa and India bear the greatest damage loss caused by rising atmospheric temperature as presented in Table~\ref{tb:table_damage_loss}. Second,  it takes regions such as China, Eurasia and Russia comparatively lower cost to reduce carbon emissions (see the estimated price for backstop technology replacing all carbon fuels in the year $2020$ in Table~\ref{tb:table_price_backstop_technology}).
	\begin{table}[htbp]
		\centering
		\small
		\setlength\tabcolsep{2pt}
		\begin{tabular}{|c|c|c|c|c|c|c|}
			\hline
			& US & EU & Japan & Russia & Eurasia & China\\ 
			\hline
			Price & 1051&1635&1635&701&701&817\\
			\hline
			& India & MidEast & Africa & LatAm & OHI & OthAsia\\
			\hline
			Price &1284&1167&1284&1518&1284&1401\\
			\hline
		\end{tabular}	
		\caption{The estimated price for backstop technology in each region that can replace all carbon fuels in the year $2020.$ }
		\label{tb:table_price_backstop_technology}	
	\end{table}

\begin{figure}[h]
			\centering
			\includegraphics[width =0.65\textwidth]{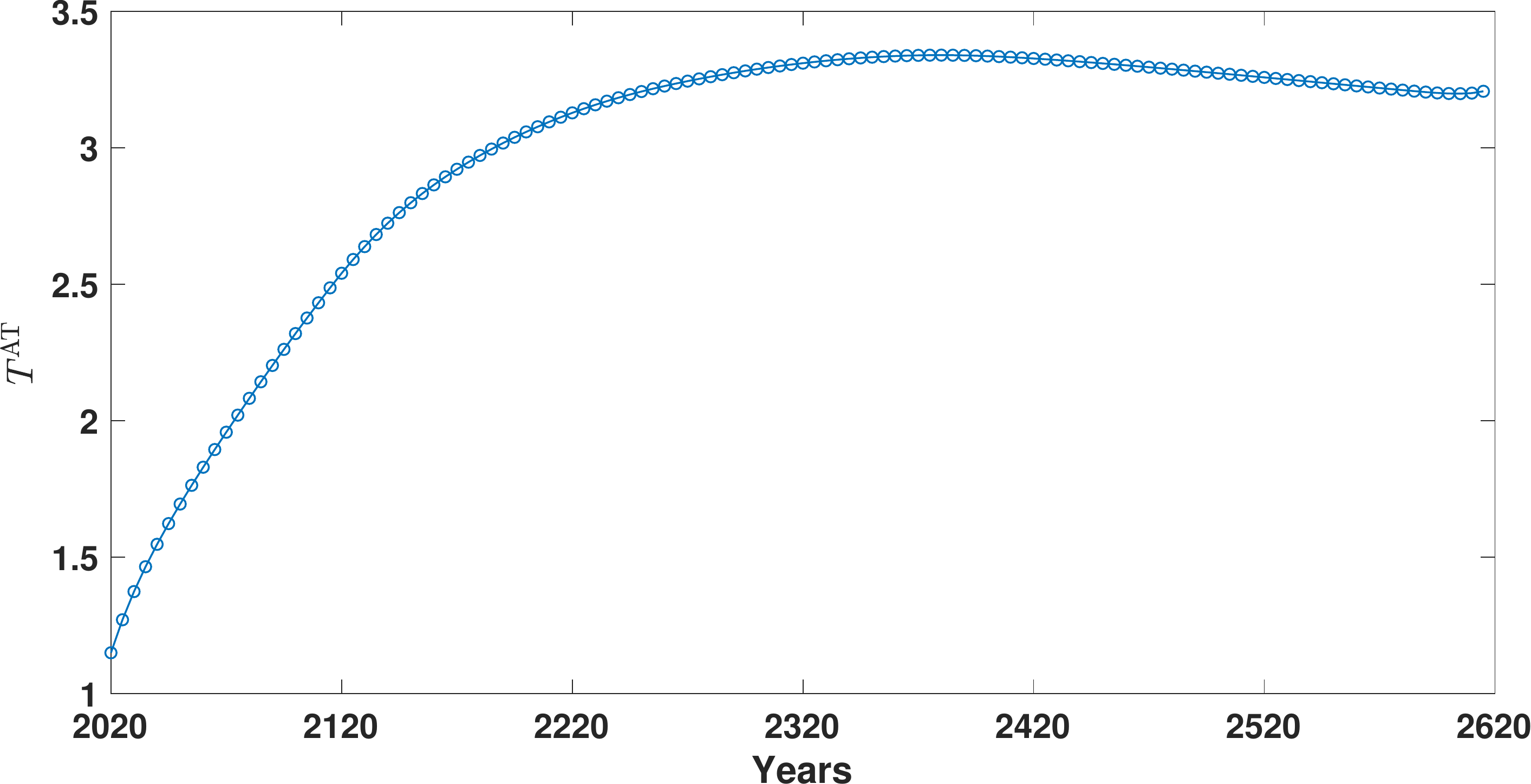}
			\caption{The atmospheric temperature deviation trajectory under $\bfU^{\rm w}$.}
			\label{fig:dg_temp}
		\end{figure}

Then, in Fig.~\ref{fig:dg_temp}, we   examine the  trajectory of the atmospheric temperature deviation when the emission-reduction rates and saving rates under $\bfU^{\rm w}$ are taken. The work of \cite{nordhaus2011estimates} also solved the OCP~\eqref{eq:swm} with the same horizon length ($600$ years)  starting from the year of $2005$. Although we are inaccessible to the exact and complete   data of results in \cite[Fig. 3]{nordhaus2011estimates}, the atmospheric temperature deviation trajectory in Fig.~\ref{fig:dg_temp} appears to be similar to that in \cite{nordhaus2011estimates}. To be specific, the atmospheric temperature deviation in \cite[Fig. 3]{nordhaus2011estimates} is around $2.8^\circ C$, whereas ours is approximately  $3^\circ C$.
		
\begin{figure}[h]
			\centering
			\includegraphics[width =0.65\textwidth]{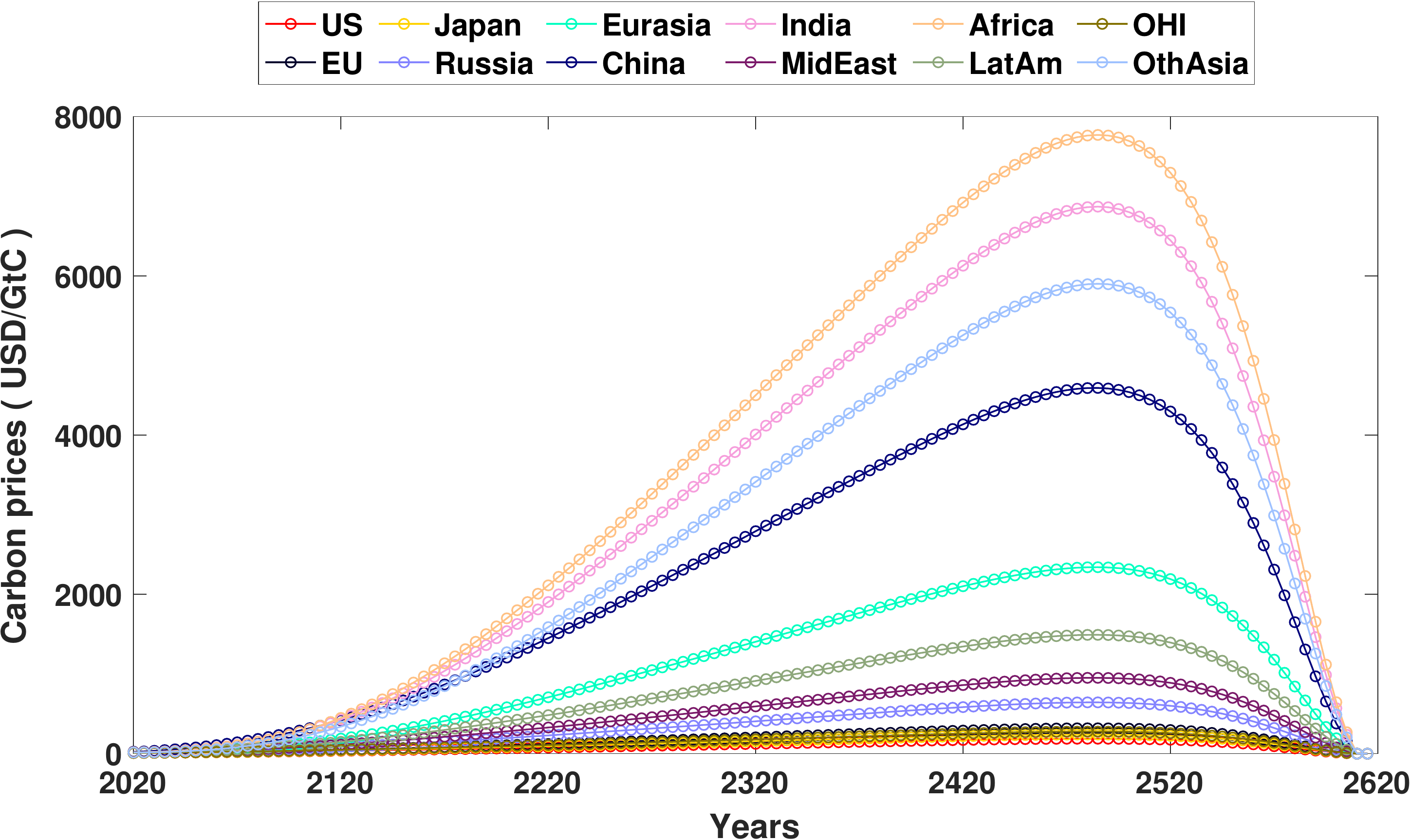}
			\caption{The social cost of ${\rm CO}_{2}$ in each region under $\bfU^{\rm w}$.}
			\label{fig:dg_scc}
		\end{figure}

Finally, in Fig.~\ref{fig:dg_scc}, the social cost of ${\rm CO}_{2}$ in each region under the global social welfare equilibrium is presented. It shows that the social cost of ${\rm CO}_{2}$ under $\bfU^{\rm w}$ in India, Africa and OthAsia is significantly higher than those in other regions. The reason for that could be that both prices for backstop technology and optimal emission-reduction rates for these regions are higher than those for other regions.
		

		\subsection{RICE Pareto Frontier}\label{subsec:Pareto_Frontier}
		
		Noting the global social welfare equilibrium is a  solution concept where the Negishi weights $c_i,i=1,\dots,n$ are calibrated in a centralized manner. Since the weight $c_i$ may greatly impact the optimal emission rates and saving rates for region $i$ and then its own social cost of carbon, the regions have intrinsic  incentives in negotiating on the distribution of the $c_i$'s. In recent years' international  climate policy forums, the divide between developed and developing regions has been one of the main barriers from a global consensus on carbon emission rates \cite{rogelj2010copenhagen,dimitrov2010inside,castro2021national}. 	In this subsection, we focus on the social welfare Pareto frontier between developed and developing regions involved in the RICE game.

\subsubsection{Solution Concept: Pareto Frontier}
		 We classify the regions in the RICE-2011 model into two clusters of regions: developed regions (US, EU, Japan and other high income countries) and developing regions (Russia, Non-Russian Eurasia, China, India, Middle East, Africa, Latin America and other Asian countries). We denote $\mcalV_{\rm developed} = \{1,2,3,11\}$ and  $\mcalV_{\rm developing} = \{4, 5, 6, 7, 8, 9, 10, 12\}$.  Correspondingly, the social welfare of the two clusters of regions is defined as, respectively, 
		$$
	\mathsf{W}_{\rm developed} = \sum_{i \in \mcalV_{\rm developed}} \mathsf{J}_i,\ \  \mathsf{W}_{\rm developing} = \sum_{i \in \mcalV_{\rm developing}} \mathsf{J}_i.
	$$ 
For the considered RICE game over these two clusters of regions, we consider the Pareto optimality.

	\begin{definition}\label{defn4}
For the RICE game under developed and developing clusters of regions, 	a decision profile  $\bfU^{\rm p}$ is a Pareto social welfare equilibrium between the developed and developing clusters   if there does not exist another decision profile  $\bfU$ such that  
		\begin{itemize}
			\item[(i)] there hold 
			\begin{align*}
			\mathsf{W}_{\rm developed} (\bfX^{\rm p},\bfU_{i}^{ \rm p},\bfU_{-i}^{\rm p})&\leq
			\mathsf{W}_{\rm developed} (\bfX,\bfU_{i},\bfU_{-i}),\\
			\mathsf{W}_{\rm developing} (\bfX^{\rm p},\bfU_{i}^{ \rm p},\bfU_{-i}^{\rm p})&\leq
			\mathsf{W}_{\rm developing} (\bfX,\bfU_{i},\bfU_{-i});
			\end{align*}
			\item[(ii)] either one of the above two inequalities holds strictly. 
		\end{itemize}		
		Here $\bfX$ and $\bfX^{\rm p}$ are the states evolved under $\bfU$ and $\bfU^{p}$, respectively. 
	\end{definition}

Based on  Lemma \ref{lemma:pareto}, we can calculate such Pareto social welfare frontier between the developed and developing clusters  by solving the family of    optimization problems for a given initial condition $\mathbf{x}_0$:
		\begin{equation}\label{eq:ocp_pareto}
			\begin{aligned}
				&  \max_{\bfU_{1},\cdots,\bfU_{n}}
				& & p \cdot \mathsf{W}_{developed} + (1-p) \cdot \mathsf{W}_{developing} \\
				&  {\rm subject \ to}
				& & \bfx(t+1) = \bff(t, \bfx(t),\bfu(t)), \ \  \bfx(0) = \bfx_{0},  t \in \mcalT \\
				&&& \bfu(t)   \in [0,1]^{24}, \ \  t \in \mcalT
			\end{aligned}
		\end{equation}
	where $p$ is selected in the interval $[0,1]$. For any fixed $p\in[0,1]$, we obtain a Pareto social welfare equilibrium, and their collection forms the Pareto frontier between the developed and developing clusters. 	 The Pareto formulation might have the potential to serve as a benchmark for the interchange of positions between developed regions and developing regions on climate policies. The parameter $p$ serves as a quantitative characterization to the allocation of responsibility for climate change mitigation between developed regions and developing regions: If $p$ is close to $0$, developed regions will take higher responsibility, as the developed regions will if $p$ is close to $1$.

\subsubsection{Results}
		We set a total of $120$ $5$-year periods. We take  $999$ linearly spaced values between $0.001$ and $0.999$ as the values of $p$. For each $p$, we obtain the Pareto solution $\bfU^{p}$ by solving the respective OCP~\eqref{eq:ocp_pareto}. We plot the social welfare Pareto frontier between developed and developing regions in Fig.~\ref{fig:Pareto_Front_Developed_Developing}. We also plot the atmospheric temperature deviation at the final time step, $T^{\rm AT}(120)$, versus the parameter $p$  in Fig.~\ref{fig:temp_deviation_vs_p}.

	\begin{figure}
	\centering
	\includegraphics[width= 0.7\textwidth]{./figs/pareto_frontier_except01-eps-converted-to.pdf}
	\caption{Social welfare Pareto frontier between developed regions and developing regions ($\times 10^{4}$ trillion USD).}
	\label{fig:Pareto_Front_Developed_Developing}
	\end{figure}
		
	\begin{figure}[h]
		\centering
		\includegraphics[width = 0.65\textwidth]{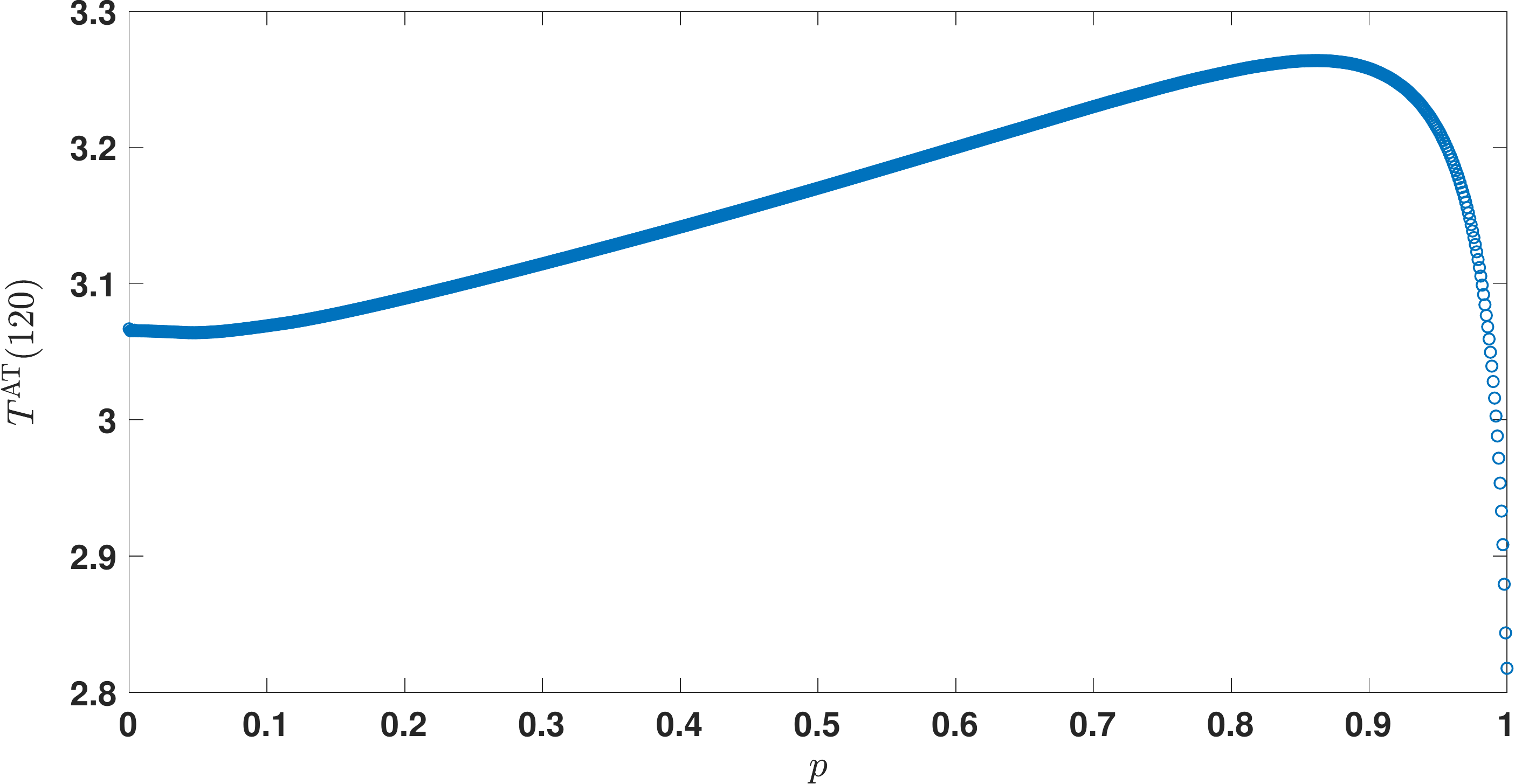}
		\caption{The atmospheric temperature deviation in the year $2620$, $T^{\rm AT}(120)$, versus the parameter $p$. }
		\label{fig:temp_deviation_vs_p}
	\end{figure}

		These results    reveal a few notable effects. 
		\begin{itemize}
			\item From Fig.~\ref{fig:Pareto_Front_Developed_Developing}, it can be seen that the values of $p$ can not drastically change the social welfare received for both the developing and developed clusters of regions. The maximal ($p=0.999$) and minimal ($p=0.001$) values  for $\mathsf{W}_{\rm developed}$ differs only by $0.35\%$; the maximal ($p=0.001$) and minimal ($p=0.999$) values  for $\mathsf{W}_{\rm developing}$ differs only by $0.56\%$. 
			\item  Fig~\ref{fig:temp_deviation_vs_p} shows that the atmospheric temperature deviation at the final time step is quite robust with respect to $p$ under the Pareto equilibrium. In fact, for all values of $p$ between $0.001$ and $0.999$,  $T^{\rm AT}(120)$, the  atmospheric temperature deviation in the year of $2620$, always falls between $2.8^\circ C$ and $3.3^\circ C$. It is also clear that a higher or a lower $p$ will both result in a lower atmospheric temperature deviation in the final year, while $p=0.862$ leads to highest terminal atmospheric temperature deviation, i.e., global warming in the worst case. 
		\end{itemize}

		\subsection{Receding Horizon  RICE}\label{subsec:MPC-RICE}
		
		The  \eqref{eq:swm}  is an OCP for a nonlinear dynamical system with non-convex cost function. When the horizon length gets larger, it becomes more difficult to solve the problem numerically. The work of \cite{kellett2019feedback} established a novel receding horizon solution to DICE, which provides robustness and computational efficiency compared to solving DICE in a long time horizon directly. In what follows, we  extend the idea of \cite{kellett2019feedback} to RICE.

		\subsubsection{Solution Concept: MPC-RICE} For the receding horizon approach, we denote the prediction horizon by $T_{rh}$ and the simulation horizon by $T_{sim}$. We introduce
		\begin{align*}
			l(t, \bfx(t),\bfu(t)) :=\sum_{i = 1}^{n} c_{i}\cdot \frac{g_{i}(C_{i}(t),L_{i}(t))}{(1+\rho_{i})^{5t}},
		\end{align*}
		 and assume a full measurement or estimate of the state $\bfx(t)$ is available at each time step $t \in \mcalT_{sim}:= \{0,1,\dots,T_{sim}\}$. The receding horizon process to approximate the global social welfare equilibrium $\bfU^{\rm w}$ is proposed in Algorithm~\ref{alg:rch}. A decision profile $\bfU^{\rm rhw} =  [\bfu^{\rm rhw}(0)^{\top},\dots,\bfu^{\rm rhw}(T_{sim})^{\top}]^{\top}$ as output of Algorithm~\ref{alg:rch} is said to be a receding horizon global social welfare equilibrium.
		
		\begin{algorithm}
			\caption{MPC-RICE}\label{alg:rch}
			{\bf{Input:}} simulation horizon $T_{sim}$; prediction horizon $T_{rh}$.   
			\begin{algorithmic}[1]
				\State $t \gets 0$
				\While{\texttt{$t \leq T_{sim}$} }
				\State	{\bf observe} $\bfx(t)$ 
				
				\State	{\bf compute} the optimal solution ${\bfu^{\ast}(s),s \in \mathcal{S}:=\{t, t+1, \dots, t+T_{rh}\}},$ to the following optimization problem over the receding horizon $\mathcal{S}$
				\begin{equation}\label{eq:rch}
					\begin{aligned}
						&\max_{ \mathbf{u}(s),\forall s \in \mcalS} \quad 	& &  \sum_{s = t}^{t+T_{rh}} l(s, \bfx(s),\bfu(s)) \\
						&  {\rm subject \ to}
						& & \bfx(s+1) = \bff(s, \bfx(s),\bfu(s)),\ \  s \in \mathcal{S},\\
						&&& \bfu(s) \in [0,1]^{24}, \ \  s \in \mathcal{S}.
					\end{aligned}
				\end{equation}
				
				\State	{\bf apply} $\bfu^{\rm rhw}(t):=\bfu^{\ast}(t)$ to RICE game
				\EndWhile
				\State \Return $\bfU^{\rm rhw}$
			\end{algorithmic}
		\end{algorithm}

		\subsubsection{Results} 
		We set the simulation horizon to be $T_{sim} =120$ (600 years), and implement MPC-RICE under prediction horizons $T_{rh} \in \{10, 20,60\}$. We also reproduce MPC-DICE with a simulation horizon of $600$ years under prediction horizons $T_{rh} \in \{10, 20,60\}$. 
		
\begin{figure}[h]
			\centering
			\includegraphics[width = 0.65\textwidth]{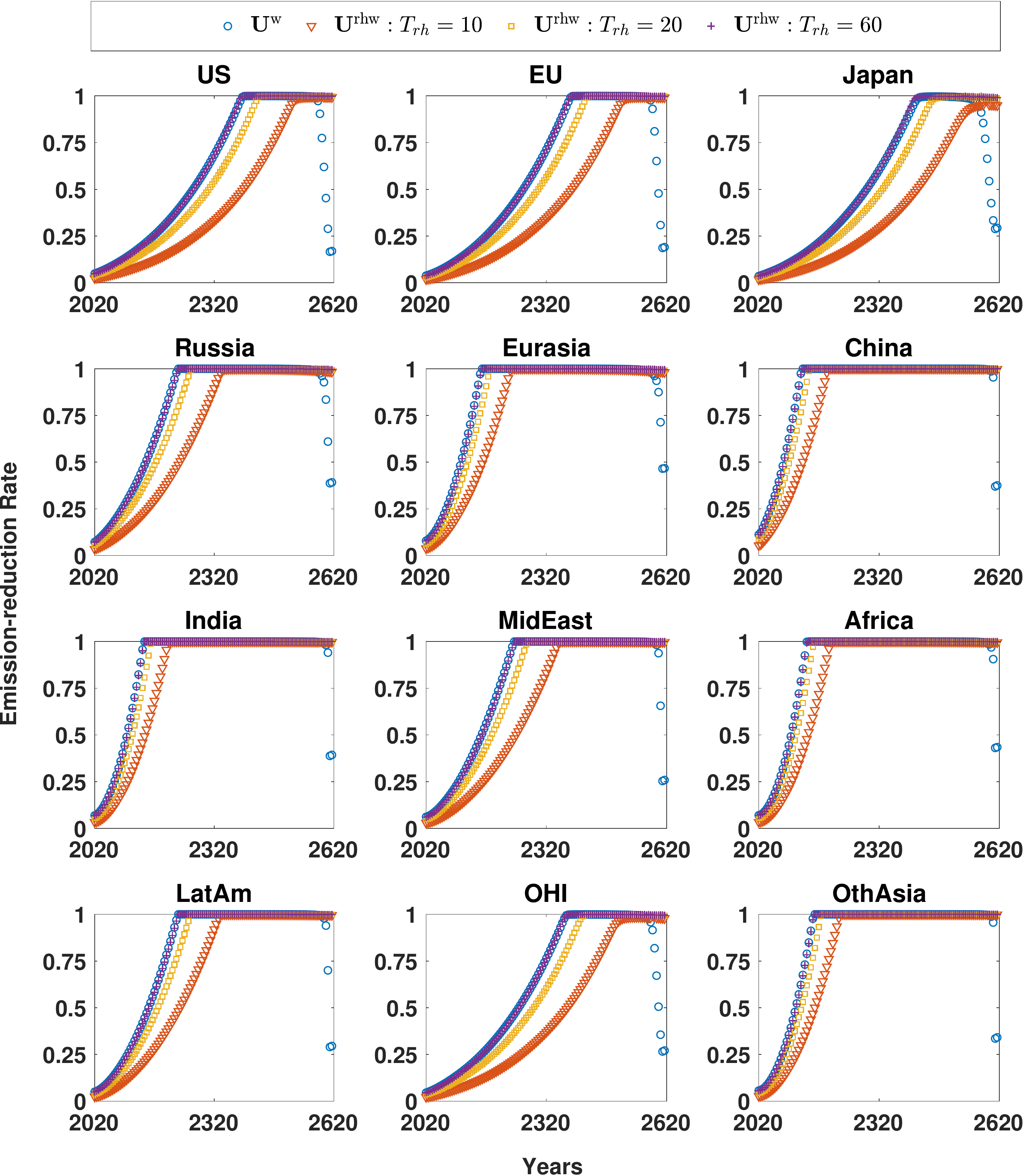}
			\caption{The comparison of each region's optimal emission-reduction rates under $\bfU^{\rm w}$ and $\bfU^{\rm rhw}$ with different prediction horizons  $T_{rh} \in \{10, 20,60\}$.}
			\label{fig:reduction_rh}
		\end{figure}

		In Fig.~\ref{fig:reduction_rh}, we plot trajectories of the optimal emission-reduction rates under global social welfare equilibrium $\bfU^{\rm w}$ and the optimal control decision $\bfU^{\rm rhw}$ obtained by MPC-RICE with $T_{rh} \in \{10, 20,60\}$. When the prediction horizon gets larger, the optimal emission-reduction rates under $\bfU^{\rm w}$ converge towards those under $\bfU^{\rm w}$ for most time steps. Moreover,  the optimal emission-reduction rates under $\bfU^{\rm rhw}$  are more steady in the sense that they do not drop back to lower levels towards the end of the simulation horizon, which is the case under $\bfU^{\rm w}$. Compared with \cite{faulwasser2018mpc}, MPC-RICE and MPC-DICE present similar approximation feature under large prediction horizons. 
		
\begin{figure}[h]
		\centering
		\includegraphics[width = 0.65\textwidth]{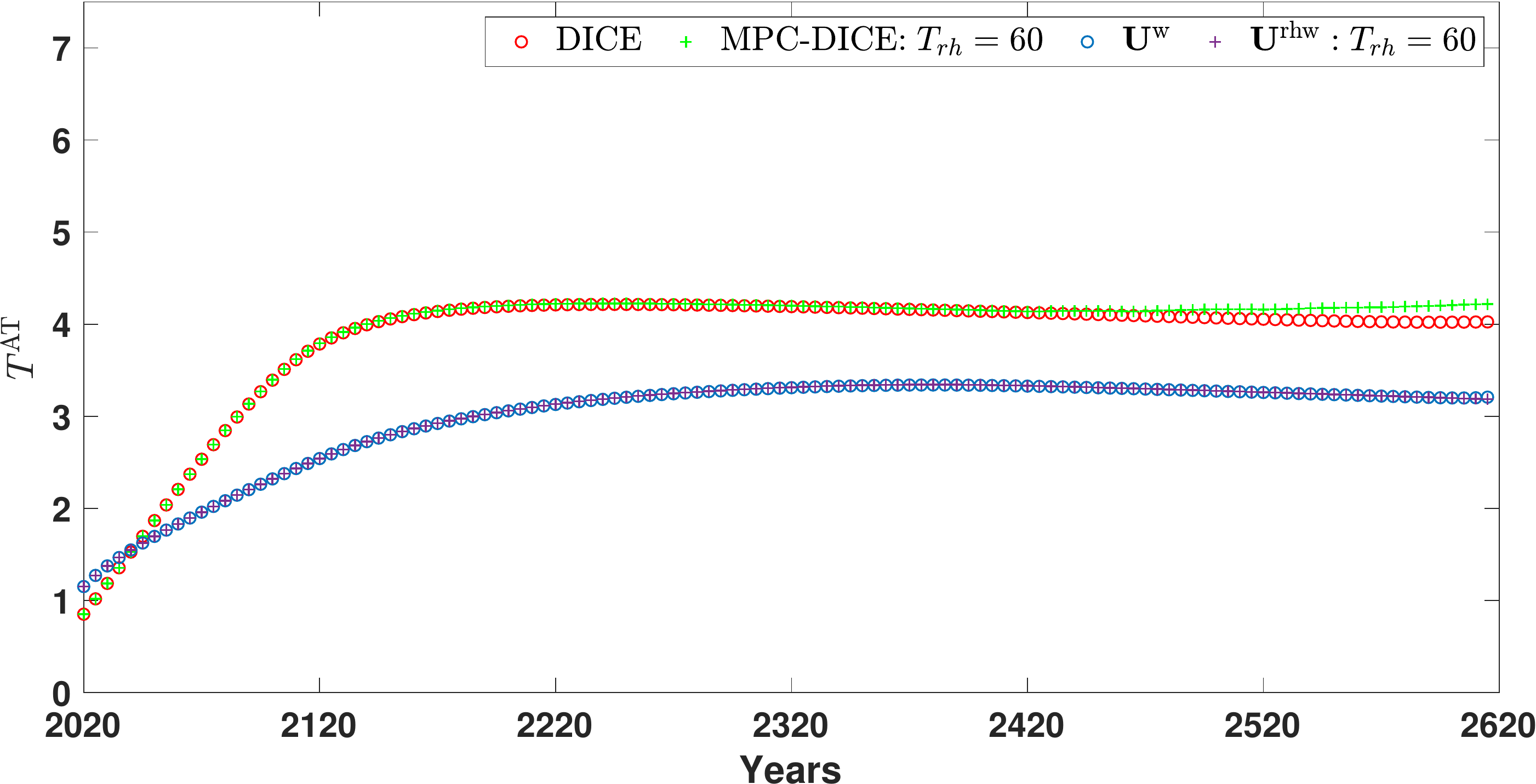}
		\caption{The comparison of atmospheric temperature deviation for the entire simulation horizon under the optimal solutions from DICE OCP, MPC-DICE with $T_{rh} = 60$, $\bfU^{\rm w}$, and $\bfU^{\rm rhw}$ with $T_{rh}  = 60$.}
		\label{fig:rh_DICE_RICE_temp}
	\end{figure}

		 In Fig.~\ref{fig:rh_DICE_RICE_temp}, we plot  trajectories of atmospheric temperature deviation for the entire simulation horizon under the optimal solutions from DICE OCP, MPC-DICE with $T_{rh} = 60$, $\bfU^{\rm w}$, and $\bfU^{\rm rhw}$ with $T_{rh}  = 60$. The trajectories of atmospheric temperature deviation under DICE  OCP and MPC-DICE with $T_{rh} = 60$ are higher than under $\bfU^{\rm w}$ and $\bfU^{\rm rhw}$ with $T_{rh} = 60$ in most time steps except at the beginning. To this extent, the RICE model and MPC-RICE raise more concerns about global warming.
		
	\section{RICE Game:  Best-response Dynamics and Open-Loop NE}\label{sec:BR-NE}
	
	In this section, we study the best-response dynamics and open-loop Nash equilibrium of the RICE game. 
	
	\subsection{Best-response Recursions  for Dynamic Games}\label{subsec:BR-DG}

	In game theory, best-response dynamics is a classical   model that describes how players strategically behave in repeated plays \cite{roughgarden2010algorithmic}. In standard best-response episodes, each player takes a best-response action for the next round that maximizes its payoff function given other players' actions in the present round of play. For certain types of games such as two-player zero-sum games and potential games,  the player actions under best-response dynamics may converge to a Nash equilibrium \cite{hofbauer2006best,christodoulou2012convergence}.

	In this subsection, we establish the  best-response recursions  for the dynamic   game introduced in Section~\ref{subsec:DTDG} with $n$ players over a    finite horizon $T$. We assume the dynamic game is repeatedly and recursively  played for $N$ episodes, where each episode consists of $T$ time steps. We thereby define the aggregated control decisions of all players in episode  $k = 1,\cdots,N $   by $\bfU^{(k)}$, and the  decisions of player $i$ in episode  $k$ by $\bfU_{i}^{(k)}$. Similarly, the decisions of players excluding player $i$ in episode  $k$ is denoted by  $\bfU_{-i}^{(k)}$. In view of the best-response for static games \cite{roughgarden2010algorithmic}, the best-response recursion of the agents in the dynamic game over the $N$ episodes are described in  the following  Recursive Best-response Algorithm for Dynamic Games (RBA-DG) as in Algorithm~\ref{alg:dbr}. The RBA-DG produces an output $\bfU^{\rm NE}_{\ast} = \bfU^{(N)}$ after $N$ episodes of updates.

	\begin{algorithm}[h]
		\caption{Recursive Best-response Algorithm for Dynamic Games (RBA-DG)}\label{alg:dbr}
		{\bf{Input:}} Episodes $N$; $c_{i}, i \in \mcalV$ .
		\begin{algorithmic}[1]
			\State  {\bf compute} an optimal cooperative solution $\bfU^{\rm c}$ by the following problem
			\begin{equation*}\label{eq:initialize}
				\begin{aligned}
					&  \max_{\bfU_{1},\cdots,\bfU_{n}}
					& & \sum_{i \in \mcalV} c_{i}\cdot J_{i}(\bfX,\bfU_{i},\bfU_{-i}) \\
					&  {\rm subject \ to}
					& & \bfx(t+1) = \bff(t, \bfx(t),\bfu(t)), \ \  \bfx(0) = \bfx_{0},  t \in \mcalT \\
					&&& \bfu(t)   \in [0,1]^{24}, \ \  t \in \mcalT. 
				\end{aligned}
			\end{equation*}
			\State {\bf let} $\bfU^{(0)}_{i} = \bfU^{\rm c}_{i}, \forall i  \in \mcalV$ 
			\State $k \gets 0$
			\While{\texttt{ k < N} }
			\For{\texttt{each player $i \in \mcalV$}}
			\State {\bf observe} $\bfU^{(k)}_{-i}$ 
			\State {\bf compute} $\bfU^{(k+1)}_{i}$  by solving the problem
			\begin{subequations}\label{eq:dbr}
				\begin{align}
					\max_{\bfU_{i}}  & \quad \mathsf{J}_i(\bfX,\bfU_{i},\bfU_{-i}),\\
					s.t. &\quad 	\bfx(t+1) = \bff(t, \bfx(t),\bfu(t)),\\
					& \quad \bfU_{-i} = \bfU^{(k)}_{-i},\\
					& \quad \bfx(0) = \bfx_{0}.
				\end{align}
			\end{subequations}
			\EndFor
			\State $\bfU^{(k+1)} = [\bfU^{(k+1)}_{1};\cdots;\bfU^{(k+1)}_{n}]$
			\State $k \gets k+1$  
			\EndWhile 
			\State \Return $\bfU^{\rm NE}_{\ast} = \bfU^{(N)}$
		\end{algorithmic}
	
	\end{algorithm}

	The implementation of RBA-DG requires two conditions. First, the initial value of the system $\bfx(0) = \bfx_{0}$ needs to be known by all players at the beginning of the process. Second, at the end of each episode $k = 1,\cdots,N$, every player should be able to  observe or know all  other players'  decision sequences over the episode $\bfU_{-i}^{(k)}$. Then the update of the player decisions for the next episode follows directly from the best-response dynamics. We present the following result, which holds true immediately from definition of open-loop Nash equilibrium.
	
	\begin{proposition}\label{prop1}
	Consider the repeatedly played $n$-player dynamic game. Let $N=\infty$ in the RBA-DG. Suppose there exists $\bfU^{\ast}$ such that there holds 
	$\lim_{k\to\infty}\bfU^{(k)} = \bfU^{\ast}.
	$
	Then $\bfU^{\ast}$ is an open-loop Nash equilibrium for the dynamic game. 
	\end{proposition}

	\subsection{Recursive Best-response for RICE}\label{subsec:BR-RICE}
	
	The RICE model describes the interplay of the various regions in the world in terms of their climate-change mitigation policy over a few centuries (for example, RICE-2011 assumes a total time horizon of $600$ years \cite[Supplementary Material]{nordhaus2011estimates}). As a result, despite the fact that the RBA-DG may be conceptually applied to RICE since RICE has been identified as a dynamic game, it is important to clarify the real-world implication of recursive best response for RICE.

	\subsubsection{Solution Concept: Regional Climate Policy Negotiations}	
		
	We propose to adopt RBA-DG over RICE as a mechanism for regional climate policy negotiations. The overall negotiations take a prescribed $N$ episodes. In each round of the negotiations, regions accept RICE as the standing model for climate-economy integration, and decide their emission reduction rates and saving rates for a fixed time horizon (e.g., 200 years). At the end of each round, all regions reveal their current planning of the emission reduction rates and saving rates for the entire time horizon to other regions. Then, during the next round of negotiations, regions get to revise their planned decisions and adopt RBA-DG as their principle of updating such planned decisions. After  the $N$ episodes of negotiations, if all regions realize none of them can unilaterally change their climate actions and gain significant increase in social welfare, such a mechanism will produce an approximate open-loop Nash equilibrium for the RICE game in view of Proposition \ref{prop1}. Such a Nash equilibrium holds higher promise of being accepted by all regions since no region is able to benefit from a revised decision when all other regions take the actions from the equilibrium.

	\subsubsection{Results}
	
	Now we implement the RBA-DG over the RICE game. We set the number of episodes  to be $N=21$. 

	\medskip

\begin{figure}[h]
	\centering
	\includegraphics[width = 0.65\textwidth]{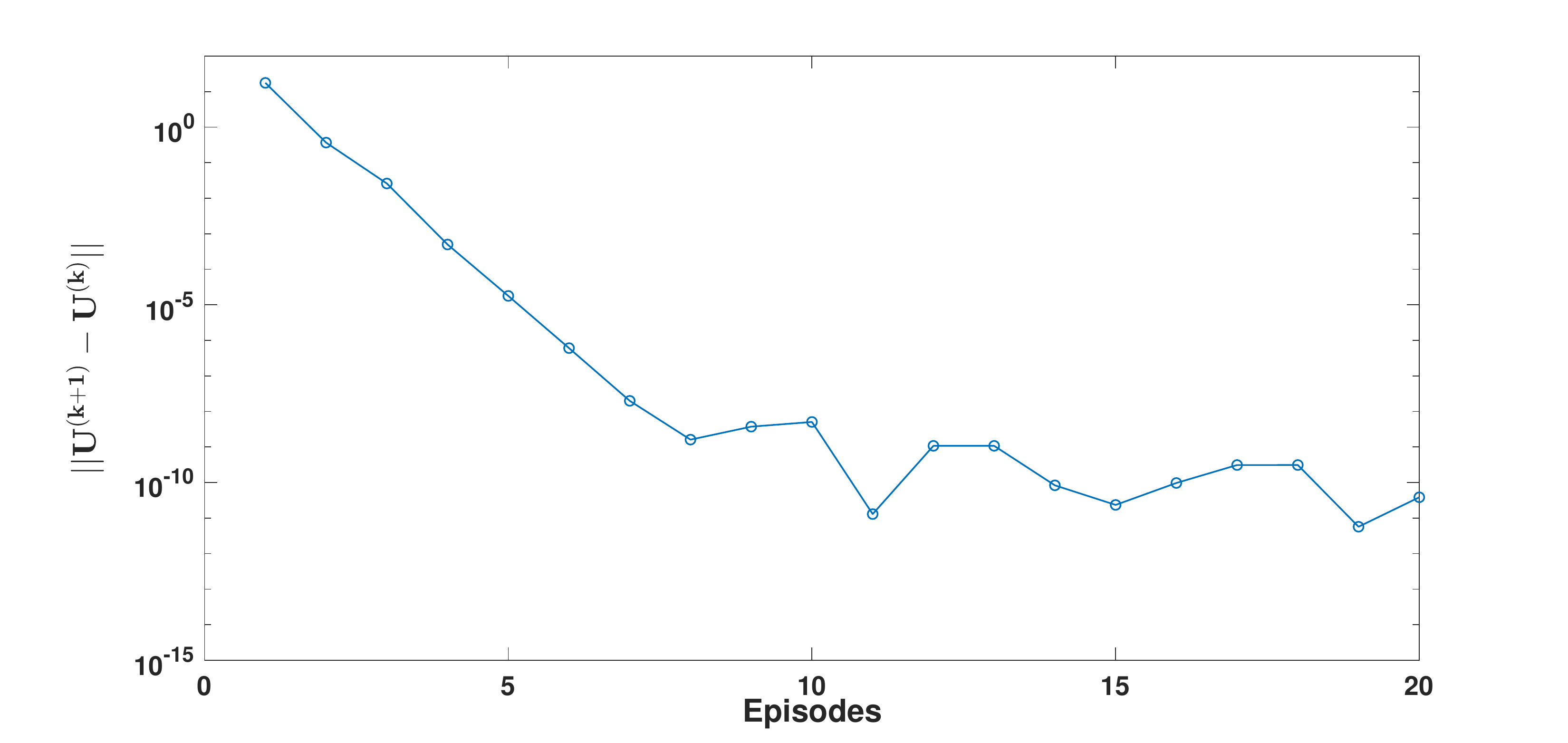}
	\caption{Convergence of $\|\bfU^{(k+1)} - \bfU^{(k)}\|$ versus episodes for RBA-DG.}.
	\label{fig:convergence_br}
	\end{figure}

	\noindent{\bf Convergence}. 	In Fig.~\ref{fig:convergence_br}, we plot the trajectory of $\|\bfU^{(k+1)} - \bfU^{(k)}\|$ versus episode $k$. The result shows that when applying RBA-DG over the RICE game, the obtained sequence of $\bfU^{(k)}$ converges to a steady point, and    after $5$ episodes, the $\|\bfU^{(k+1)}- \bfU^{(k)}\|$ has become very close to zero. From Proposition \ref{prop1}, this implies that the RBA-DG may also serve as an efficient algorithm for computing the open-loop Nash equilibrium of the RICE game. 
		
	\medskip

\begin{figure}[h]
		\centering
		\includegraphics[width = 0.65\textwidth]{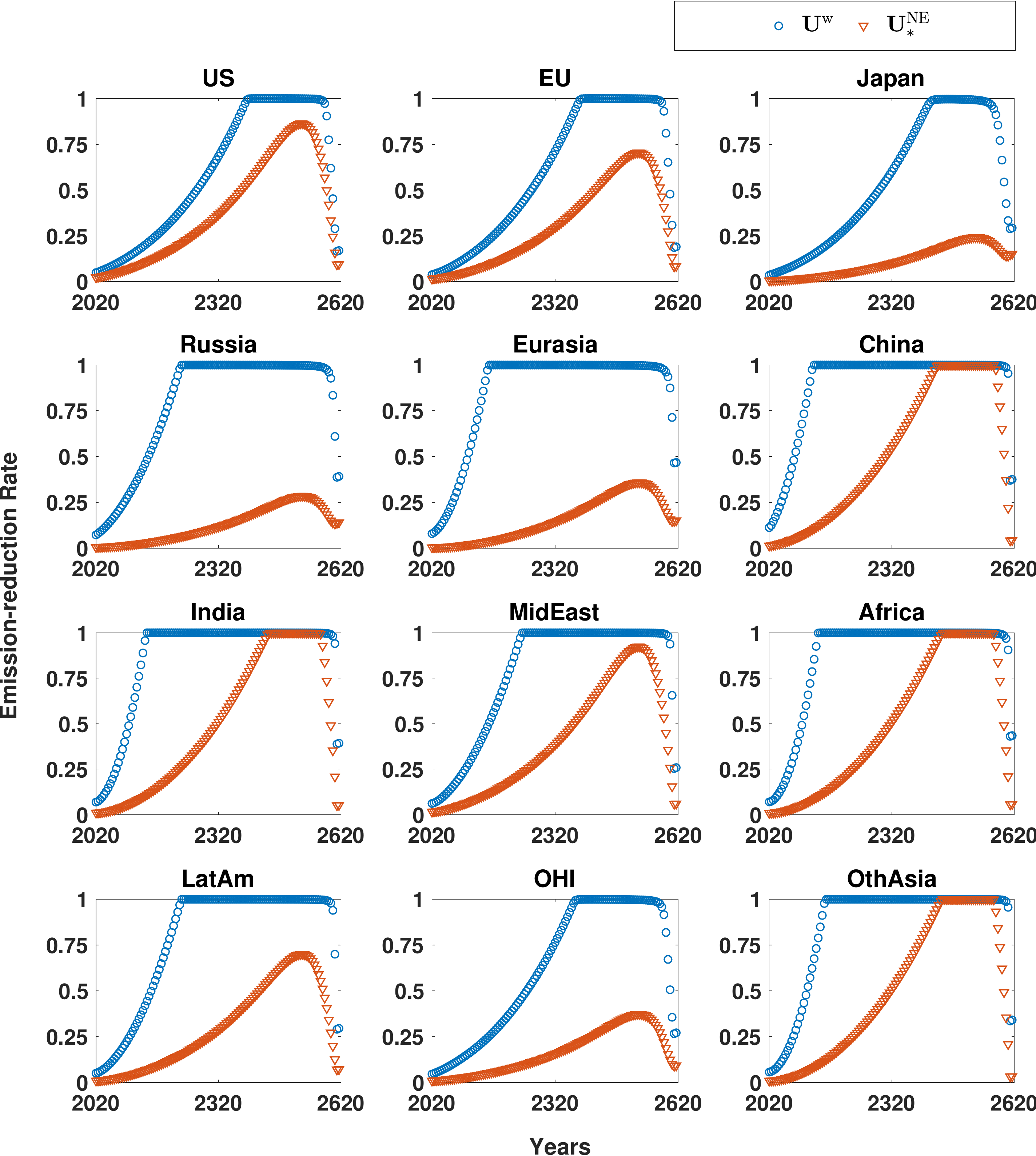}
		\caption{The comparison of each region's optimal emission-reduction rates under $\bfU^{\rm w}$ and $\bfU^{\rm NE}_{\ast}$.}
		\label{fig:reduction_br}
	\end{figure}

	\noindent{\bf Cooperation vs Competition}. In Fig.~\ref{fig:reduction_br}, we plot the comparison of each region's optimal emission-reduction rates under global social welfare equilibrium $\bfU^{\rm w}$ (cooperative solution) and the optimal control decision $\bfU^{\rm NE}_{\ast}$ (competitive solution) obtained  by RBA-DG (Algorithm~\ref{alg:dbr}). From Fig.~\ref{fig:reduction_br}, it is not surprising that each region's optimal emission-reduction rates under $\bfU^{\rm NE}_{\ast}$ are significantly lower than those under $\bfU^{\rm w}$. This is a direct reflection that under Nash equilibrium of the RICE game, the regions are working towards maximizing their own social welfare instead of a collective social welfare of all regions. As a result,  each region has   incentives to reduce its emission-reduction rate and its emission abatement cost, and thus improve its social welfare. These results partially show the rationale behind regions leaving signed  climate treaties, e.g., Canada opt-out from Kyoto Protocol \cite{nwanguma2016analysis} in $2012$, and   the United States formally quit Paris Agreement \cite{rhodes2017us} in $2020$ because a less aggressive regional policy leads to higher economic benefit even facing climate change damages.
	
\begin{figure}[h]
	\centering
	\includegraphics[width = 0.65\textwidth]{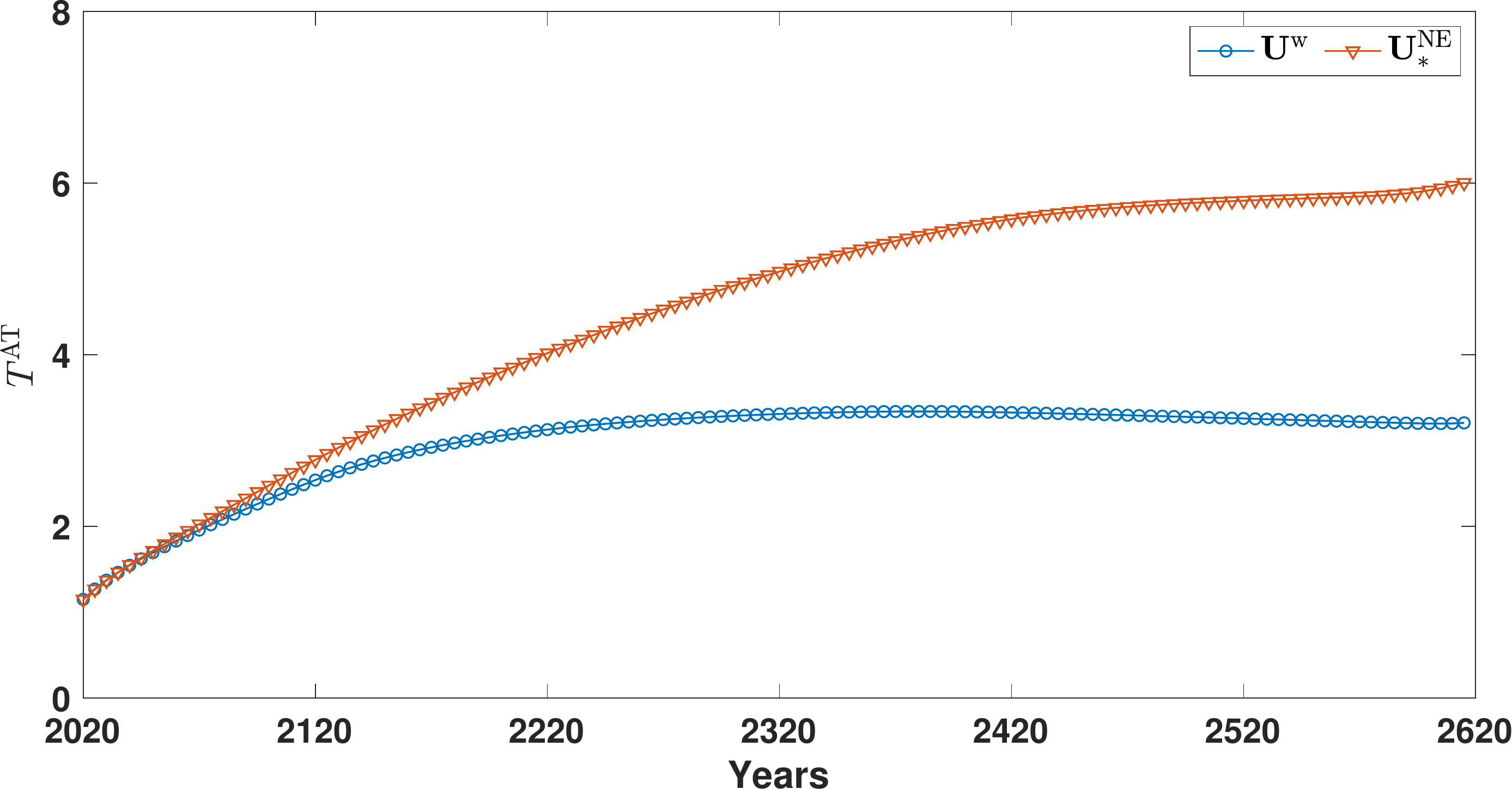}
	\caption{The comparison of the atmospheric temperature deviation trajectories under $\bfU^{\rm w}$ and $\bfU^{\rm NE}_{\ast}$.}
	\label{fig:temp_br}
	\end{figure}

	In Fig.~\ref{fig:temp_br}, we plot the comparison of the atmospheric temperature deviation trajectories under $\bfU^{\rm w}$ and $\bfU^{(\rm NE)}_{\ast}$. From the results, with competition, substantially higher atmospheric temperature deviation occurs under $\bfU^{(\rm NE)}_{\ast}$, as a consequence of lower emission-reduction rates. For example, in $2620$, the atmospheric temperature deviation resulted from  RBA-DG  is around $6^\circ$C while the atmospheric temperature deviation under global social welfare  equilibrium is about $3^\circ$C.
	
	\section{RICE Game: Receding Horizon Feedback Decisions}\label{sec:RHFD}

	In this section, we present a framework for dynamic games where in a single play over the time horizon,  players observe the underlying dynamic process and other players' actions, and then apply a receding-horizon feedback decision making   for a prediction horizon. Then we apply this receding horizon feedback process to RICE game in order to capture the competitive nature of regional climate policies, and show how such competitions have impact on the global climate dynamics. 
 
	\subsection{Receding Horizon Feedback Decisions for Dynamic Games} \label{subsec:RHFD-DG}

	Consider  the dynamic   game introduced in Section~\ref{subsec:DTDG} with $n$ players over a    finite horizon $T$. The game is played only once, and the players take the following feedback decision process in the receding horizon sense. 

	\medskip

	\noindent{\bf The Receding Horizon Feedback Decision Process.} The players apply a receding horizon approach and compute their feedback decisions $\mathbf{u}_i(t)$. At each time $t=0,\dots,T-1$, each player $i$ observes other players'  played action   $\mathbf{u}_{-i}(t)$ and the system state $\mathbf{x}(t)$. Then, every player $i$ assumes that $\mathbf{u}_{-i}(t)$ will continue to be played over $[t+1,t+T_{rh}]$, and therefore decides its best feedback decision plan $\mathbf{u}_{i|t+1\to t+T_{rh}}^{\rm RHP} (\mathbf{x}(t),\mathbf{u}_{-i}(t))$, where 
	$\mathbf{u}_{i|t+1\to t+T_{rh}}^{\rm RHP}$ maximizes the cumulative payoff of player $i$ over the time horizon  $[t+1,t+T_{rh}]$ conditioned on that $\mathbf{u}_{-i}(t)$  be played over $[t+1,t+T_{rh}]$. Finally, each player $i$ plays the first planned decision $\mathbf{u}_{i|t+1}^{\rm RHP}$ in $\mathbf{u}_{i|t+1\to t+T_{rh}}^{\rm RHP}$ for the step $t+1$, and the process moves forward recursively.  Denoting $\mathbf{u}_i^{\rm RHF}(t)$ as the actions  generated by the receding horizon feedback decision process, clearly there is an underlying feedback law $\pi_i$ such that
	$$
	\mathbf{u}_i^{\rm RHF}(t)=\pi_i(t,\mathbf{x}(t),\mathbf{u}_{-i}^{\rm RHF}(t)).
	$$
	Note that the decisions $\bfu_{i}^{\rm RHF}(t)$ are actually played by the players at each time $t$. 
	The resulting  collective player decisions over the entire time horizon is written as $\mathbf{U}^{\rm RHF}$. The computational process of this receding horizon feedback decision framework is presented in the following  Receding Horizon Feedback Algorithm for Dynamic Games (RHFA-DG). 

	\begin{algorithm}[h]
		\caption{Receding Horizon Feedback Algorithm for Dynamic Games (RHFA-DG)}\label{alg:drh}
		
		{\bf{Input:}} simulation horizon $T_{sim}$; prediction horizon $T_{rh}$; $c_{i}, i \in \mcalV$. 
		
		\begin{algorithmic}[1]
			\State {\bf compute} an optimal cooperative solution $\bfU^{\rm c}$ by the following problem
			\begin{equation*}
				\begin{aligned}
					&  \max_{\bfU_{1},\cdots,\bfU_{n}}
					& & \sum_{i \in \mcalV} c_{i}\cdot J_{i}(\bfX,\bfU_{i},\bfU_{-i}) \\
					&  {\rm subject \ to}
					& & \bfx(t+1) = \bff(t, \bfx(t),\bfu(t)), \ \  \bfx(0) = \bfx_{0},  t \in \mcalT \\
					&&& \bfu(t)   \in [0,1]^{24}, \ \  t \in \mcalT. 
				\end{aligned}
			\end{equation*}
			\State {\bf let} $\bfU^{\rm RHF}_{i}(0) = \bfU^{\rm c}_{i}(0), \forall i  \in \mcalV$ 
			
			\State $t \gets 0$
			
			\While{\texttt{$t \leq T_{sim}$} }
			\For{\texttt{each player $i \in \mcalV$}} 
			\State 	{\bf Take action}  $\bfU^{\rm RHF}_{i}(t)$
			\EndFor
			
			\For{\texttt{each player $i \in \mcalV$}} 
			\State {\bf observe} $\bfx(t)$ and $\bfU^{\rm RHF}(t)$
			
			\State {\bf compute} $\bfx(t+1)$ according to~\eqref{eq:DTDG_dynamics} 
			
			\State {\bf assume} all players $j \in \mcalV/\{i\}$ will continue to play $\bfU^{\rm RHF}_{-i}(t)$ over $[t+1,t+T_{rh}]$
			
			\State {\bf compute} its optimal solution $\bfu_{i|t+1\to t+T_{rh}}^{\rm RHP}$  to the following receding horizon optimization problem
			\begin{subequations}\label{eq:drh}
				\begin{align}
					\max_{ \bfu_{i|t+1\to t+T_{rh}}} \quad & \sum_{s = t+1}^{t+T_{rh}} g_{i}(s, \bfx(s),\bfu_{i}(s),\bfu_{-i}(s)) \\
					s.t. \quad & \bfx(s+1) = \bff(s, \bfx(s),\bfu(s)),\\
					& \bfu_{-i}(s) = \bfU^{\rm RHF}_{-i}(t),  \ \ s \in [t+1,t+T_{rh}].
				\end{align}
			\end{subequations}

			\State {\bf plan} $\bfU^{\rm RHF}_{i}(t+1) = \bfu_{i|t+1}^{\rm RHP}$
			\EndFor
			\EndWhile
		\end{algorithmic}
	\end{algorithm}
	
	\subsection{Receding Horizon Feedback for RICE}\label{subsec:RHF-RICE}
	
	In the current climate-change mitigation measures, there has been no international  consensus on the emission reduction rates and the saving rates for different  regions or for the globe collectively, despite the successful adoption of global or regional  climate-change agreements such  as the Paris Agreement.  In fact, most climate treaties are neither substantial nor mandatory \cite{cass2012failures}. The objective  of such treaties has  not been the actual  emission-reduction rate, but rather to establish environmental norms at the international level. The hope is that such international environmental norms may then be translated into the domestic climate policies according to each region's political processes \cite{cass2016measuring}. 
	
	In the real world, regions revise their climate-change policies from time to time and attempt to compete with each other while acknowledging the importance of climate-change mitigation. Indeed, such revisions of regional climate-change policies may depend on many factors such as the public opinion shifts, government changes, and new international environmental norms. In the end, the nature of regional competitions persists   as aggressive climate-change mitigation policy for a region may benefit other regions economically in the short run, although collectively aggressive  climate-change mitigation policy benefits every region as shown from the comparison between  cooperative social welfare equilibrium and the approximate open-loop Nash equilibrium. Here, we propose to apply RHFA-DG on  RICE as an attempt to model the real-world regional climate-policy evolution given the economically competitive nature of such policies.

	\subsection{Results}
	We implement the RHFA-DG over the RICE game. We set the simulation horizon and prediction horizon to be $T_{sim} = 120$ and $T_{rh} = \{5,10,20\}$.

\begin{figure}[h]
	\centering
	\includegraphics[width =  0.65\textwidth]{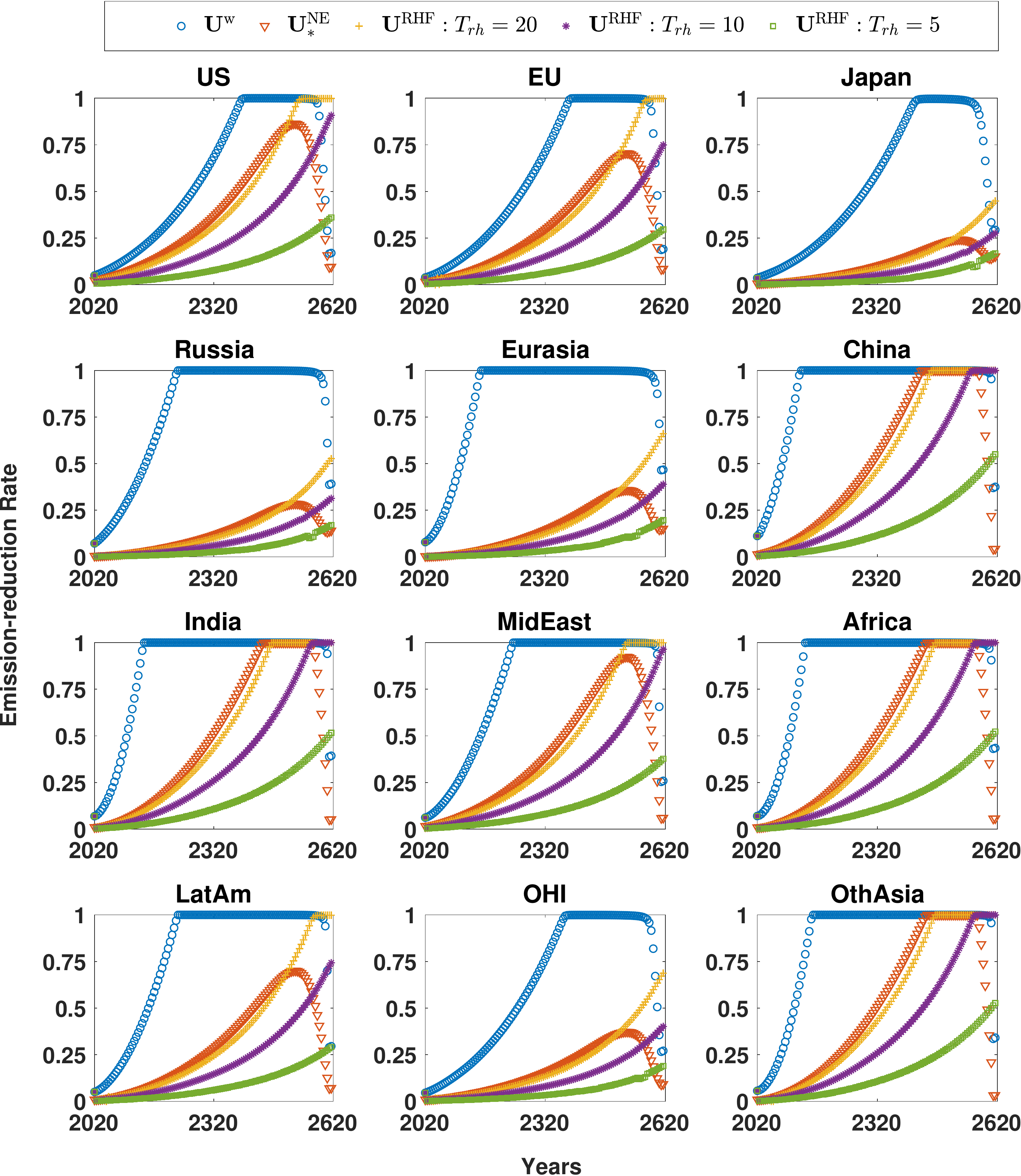}
	\caption{The comparison of each region's optimal emission-reduction rates under $\bfU^{\rm w}$, $\bfU^{\rm NE}_{\ast}$, and $\bfU^{\rm RHF}$ with $T_{rh} = \{5,10,20\}$.}
	\label{fig:reduction_drh}
	\end{figure}

	In Fig.~\ref{fig:reduction_drh}, we plot the comparison of each region's optimal emission-reduction rates under global social welfare equilibrium $\bfU^{\rm w}$, the approximated open-loop Nash equilibrium $\bfU^{\rm NE}_{\ast}$ solved from RBA-DG (Algorithm~\ref{alg:dbr}), and the receding horizon feedback decisions $\bfU^{\rm RHF}$ solved from RHFA-DG (Algorithm~\ref{alg:drh}) with $T_{rh} = \{5,10,20\}$. First, with competition, the receding horizon feedback decisions of emission-reduction rates solved from RHFA-DG are significantly lower than those under global social welfare equilibrium. Second, due to receding horizon and myopic assumption of other regions' future decisions, the emission-reduction rates from RHFA-DG are lower than  those obtained from RBA-DG in the early time steps. However, interestingly, the former continues to climb in the final time steps, while the latter drops back to lower levels in the end of simulation horizons.
	
\begin{figure}[h]
	\centering
	\includegraphics[width = 0.65\textwidth]{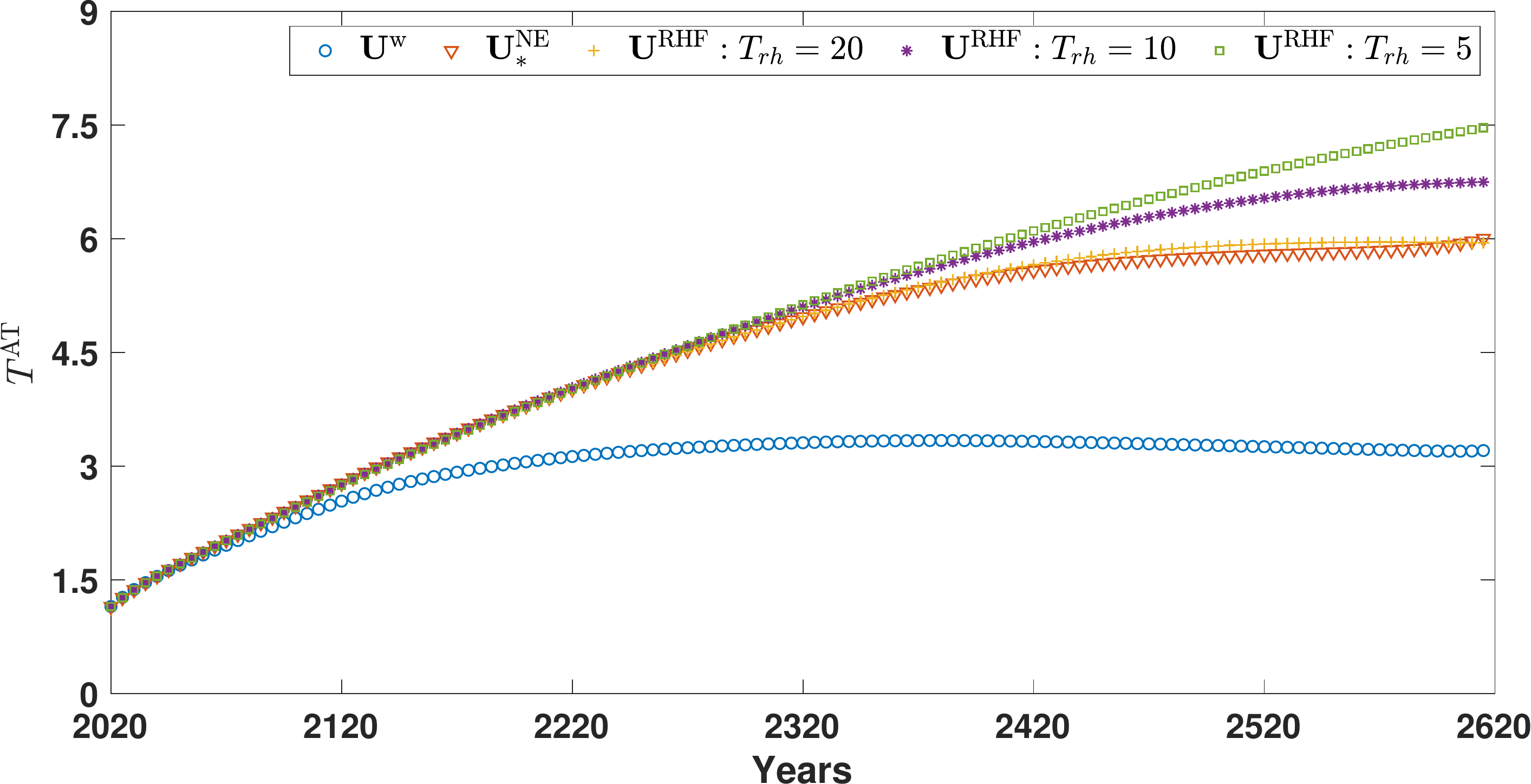}
	\caption{The comparison of the atmospheric temperature deviation trajectories under $\bfU^{\rm w}$, $\bfU^{\rm NE}_{\ast}$, and $\bfU^{\rm RHF}$ with $T_{rh} = \{5,10,20\}$.}
	\label{fig:temp_drh}
	\end{figure}

	In Fig.~\ref{fig:temp_drh}, we plot the comparison of the atmospheric temperature deviation trajectories under $\bfU^{\rm w}$, $\bfU^{\rm NE}_{\ast}$, and $\bfU^{\rm RHF}$ with $T_{rh} = \{5,10,20\}$. There are several observations on the results. First, the atmospheric temperature deviation under $\bfU^{\rm RHF}$ is substantially higher than that under $\bfU^{\rm w}$, and slightly higher than that under $\bfU^{\rm NE}_{\ast}$. Second, for $\bfU^{\rm RHF}$ with $T_{rh} = \{5,10,20\}$, the smaller the prediction horizon is, the higher the atmospheric temperature deviation is. This implies that longer prediction horizon forces the regions to take the long-term climate damages more into their receding horizon decisions, and therefore leads to better climate-change mitigation.

	\section{Conclusions}\label{sec:conclusions}
	In this paper, we investigated how cooperation and competition arise in regional climate policies under the RICE framework from the perspective of game theory and optimal control. The results revealed how game theory may be used to facilitate international negotiations towards consensus on regional climate-change mitigation policies, as well as how cooperative and competitive regional relations shape climate change for our future.
	
	\bibliographystyle{plain}
	\bibliography{ref}
	
\end{document}